%% file: bratteli10_C14.tex
\documentclass[11pt,a4paper]{article}

\include{macros}

\title{Tiling groupoids and Bratteli diagrams II: \protect\\ structure of the orbit equivalence relation}
\author{ A. Julien, J. Savinien\\
{\small Institut Camille Jordan, Universit\'e Lyon I, France.}
}

\begin{document}

\maketitle

\begin{abstract}
In this second paper, we study the case of substitution tilings of $\RM^d$.
The substitution on tiles induces substitutions on the faces of the tiles of
all dimensions \(j=0, \ldots, d-1\).
We reconstruct the tiling's equivalence relation in a purely combinatorial way
using the $AF$-relations given by the lower dimensional substitutions.

We define a Bratteli multi-diagram $\Bb$ which is made of the Bratteli diagrams
$\Bb^j, j=0, \ldots d$, of all those substitutions.
The set of infinite paths in $\Bb^d$ is identified with the canonical
transversal $\Xi$ of the tiling.
Any such path has a ``border'', which is a set of tails in $\Bb^j$ for some $j\le d$,
and this corresponds to a natural notion of border for its associated tiling.
We define an {\em \'etale} equivalence relation $\rel_\Bb$ on $\Bb$  by saying that two
infinite paths are equivalent if they have borders which are tail equivalent in
$\Bb^j$ for some $j\le d$.
We show that $\rel_\Bb$ is homeomorphic to the tiling's equivalence relation
$\rel_\Xi$.
\end{abstract}


\tableofcontents



\section{Introduction}

\subsection{Context}
In this article, we present a generalized version of Bratteli diagrams and use it to
encode, in a purely combinatorial way, the orbit equivalence relation 
on the transversal of substitution tiling spaces.
The structure of the diagram also brings an understanding of the structure
of the equivalence relation.

\subparagraph{Bratteli diagrams}
Bratteli diagrams have been efficiently used to encode $\ZM$-actions on the
Cantor set.
A Bratteli diagram $\Bb$ is given by sets of edges and vertices:
\[
 \Ee = \bigcup_{n \in \NM}{\Ee_n} \ ; \quad \Vv = \bigcup_{n \in \NM}{\Vv_n},
\]
and the set of infinite paths on the diagram, $\Pi_\infty$, is a closed subset of
$\prod_{n \in \NM}{\Ee_n}$ (a path being a sequence of \emph{composable} edges).
\begin{figure}[h!]
\[
\xymatrix{
&   &  b \ar@{-}[rrdd] & &  b \ar@{-}[rrdd]   & &  b \ar@{-}[rrdd]  & &  b \ar@{.}[dr] & \\
&\circ \ar@{-}[ur] \ar@{-}[dr] & & & & & & & &  \\
&   &  a \ar@{-}[rr]  \ar@{-}[uurr] & &  a\ar@{-}[rr] \ar@{-}[uurr]  & &  a \ar@{-}[rr] \ar@{-}[uurr] & &  a \ar@{.}[r] \ar@{.}[ur] & \\
}
\]
\caption{\small{A self-similar Bratteli diagram (root on the left): for all $n\ge 1$, \(\vtx_n \cong \{ a,b\}\), and \(\edg_n\) corresponds to the substitution \(a \rightarrow ab, b \rightarrow a\) (see Definition~\ref{bratteli10.def-substitution}}}
\label{bratteli10.fig-BdiagFibo}
\end{figure}
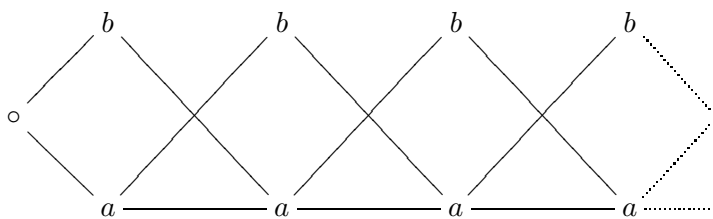

Under some conditions on the diagram, $\Pi_\infty$ is a Cantor set.
An additional structure on $\Bb$ gives a partial order on the set of paths.
With respect to this order, the ``successor'' function is well-defined, and defines a minimal action via the Vershik map \cite{VL92}.
Conversely, to any minimal action on the Cantor set, one can associate a Bratteli
diagram with order, such that the actions are conjugate.
Therefore, ordered Bratteli diagrams provide  \emph{combinatorial models} for
minimal $\ZM$-actions on the Cantor set \cite{HPS92,GPS95,DHS99}.

\subparagraph{Tilings}
A particular case of minimal $\ZM$-action on a Cantor set is the following.
Consider a bi-infinite word $w \in \proto^\ZM$, with $\proto$ a finite set
of symbols (alphabet).
Then, consider all translates of $w$ (\emph{i.e.} its orbit by the shift),
and take a closure in $\proto^\ZM$:
\[
 \Xi_w := \overline{\{ \sigma^n (w) \ ; \ n \in \ZM \}}.
\]
Under suitable conditions on $w$, the $\ZM$-dynamical system $(\Xi_w, \sigma)$ is
minimal, and $\Xi_w$ is a Cantor set.
In the case where $w$ is a word obtained by a symbolic \emph{substitution} (see Section~\ref{bratteli10.sect-stilings}), the associated Bratteli diagram can be chosen to be  self-similar.
Figure~\ref{bratteli10.fig-BdiagFibo} shows the example of the construction of Definition~\ref{bratteli10.def-substitution} for Fibonacci substitution: \(\Aa = \{ a, b \}\) and \( \omega(a) = ab, \omega(b) = a\).
\begin{figure}[htp]
 \begin{center}
\includegraphics{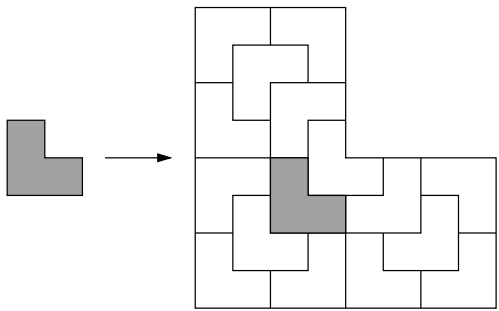}
   \caption{{\small A process of inflation and substitution (chair tiling). A whole
tiling can be obtained as a fixed point of this map.}}
   \label{bratteli10.fig-subst}
 \end{center}
\end{figure}

Tilings, and in particular substitution tilings (see Figure~\ref{bratteli10.fig-subst} for an example), are higher dimensional analogues. 
Given a tiling of $\RM^d$, let $\Omega$ be its tiling space (a closure of its family of translates, see Section~\ref{bratteli10.ssect-stilings}).
It is an $\RM^d$-dynamical system which is minimal under some conditions, and one can chose a transversal $\Xi$ to the $\RM^d$ flow (Definition~\ref{bratteli10.def-Xi}).
If $d=1$, there is a first-return map on $\Xi$ which implements a $\ZM$-action, and the case is similar to the case described above.
In higher dimension, there is no longer a group action, but a groupoid replaces $\ZM$.
This is the groupoid of the \emph{orbit equivalence relation} \(\rel_\Xi \subset \Xi \times \Xi\) defined by:
\[
(T, T') \in \rel_\Xi \quad \iff \quad T' = T + a \,, \quad \textrm{for some } a\in \RM^d\,,
\]
and with a certain topology (called {\it \'etale}, and which is not the product topology, see Definition~\ref{bratteli10.def-ERhull}).
Therefore, in dimension $d \geq 2$, the problems need then to be rephrased in terms of
equivalence relations.

The questions we address in this series of papers are the following.
\begin{itemize}
\item  Given a tiling space transversal, is it possible to give a
combinatorial description of its orbit equivalence relation in terms
of Bratteli diagrams?
\item Is it possible to describe precisely the structure of this equivalence
relation?
\end{itemize}

\subsection{Previous work}
In a previous paper~\cite{BJS10}, we addressed these questions for general tilings.
It is a known construction~\cite{Kel95,For97} that one can associate a Bratteli diagram to any given tiling space transversal $\Xi$.
On this diagram, if two paths of $\Pi_\infty$ eventually agree, then the corresponding tilings are in the same orbit.
This gives a \emph{strict} sub-equivalence relation of $\rel_\Xi$, called the \emph{tail} or \emph{$AF$-equivalence relation}, and denoted $\Rr_{AF}$, :
\[
  \Rr_{AF} \subsetneq \Rr_\Xi.
\]
In~\cite{BJS10}, we added ``horizontal edges'' to the diagram, as well as labels
on the edges. Using these data, it was possible to reconstruct the whole
equivalence relation: to add the missing parts to $\Rr_{AF}$ in order to recover
$\Rr_\Xi$.
This diagram allows to recover $\Rr_\Xi$ by defining a generalized tail equivalence
relation. However, it fails to be purely combinatorial, as the labels carry
geometric information (they are essentially translation vectors).

\subsection{Present work}
In this paper, we give a construction which holds \emph{a priori}  for substitution tilings (see Section~\ref{bratteli10.ssect-perspectives} for a discussion), but gives a much better understanding of what the ``missing parts'' are.
Furthermore, the construction presented here is purely combinatorial.
In this sense, this work generalizes the combinatorial representation of minimal $\ZM$-actions on the Cantor set by ordered Bratteli diagrams.

We assume that tiles have a good notion of faces in all dimension: see Hypothesis~\ref{bratteli10.hyp-CW} and Section~\ref{bratteli10.ssect-indsub}, as well as Figure~\ref{bratteli10.fig-faces}.
The diagram we use here is a multi-diagram built as follows.
Given the transversal $\Xi$ of a tiling space of dimension $d$ associated with a
substitution $\omega_d$,  build its usual Bratteli diagram $\mbrat{d}$.
Let then $\omega_j$ be the substitution induced by $\omega_d$ on the $j$-dimensional faces of the tiles.
For all $0 \leq j \leq d-1$, build $\mbrat{j}$ the Bratteli diagram of $\omega_j$.
All these diagrams are then linked by {\em horizontal edges}, which encode adjacencies (how a face of dimension $j$ contains sub-faces of dimension $j-1$), see Section~\ref{bratteli10.ssect-hor} and Figure~\ref{bratteli10.fig-hor}.

As before, there is a homeomorphism between the set $\Pi_\infty^d$ of infinite paths in $\mbrat{d}$, and the transversal.
We define {\em borders of a path} $x\in\Pi_\infty^d$ (Definition~\ref{bratteli10.def-borderpath}), as sets of tails in the diagrams of lower dimensions (infinite paths that need not start at depth $1$) -- the rule on how to derive them being given by horizontal edges.
We denote by $\bo^j(x)$, $j\le d$, the set of tails in $\Bb^j$ derived from $x$, and call it its {\em $j$-th border}.
The smallest $j$ for which $\bo^j(x)$ is non empty, is called the {\em border dimension} of $x$, and written \(\bd(x)\).
This is equivalent to the following natural notion of border for tilings (Proposition~\ref{bratteli10.prop-bordertiling2}):
for $T\in\Xi$ define its $j$-th border (Definition~\ref{bratteli10.def-bordertiling}) as
\[
\bo^j(T) = \bigcap_{n\in\NM} \lambda^n \; \omega_d^{-n}(T)^j \,,
\]
where \(T^j\) denote the union of the $j$-faces of the tiles of $T$, and $\lambda$ is the dilation factor of the substitution; this is discussed in details in Section~\ref{bratteli10.ssect-border}.

If \(\bd(x)=j\), then any two tails in $\bo^j(x)$  are tail-equivalent in $\Bb^j$ (Lemma~\ref{bratteli10.lem-bdim}).
We can then define an equivalence relation $\rel_\Bb$ on $\Pi_\infty^d$ as follows.
We say that two paths \(x,y \in \Pi_\infty^d\) are {\em border equivalent}, and write \(x\sim y\), if they have the same border dimension \(j \), and their borders are tail equivalent in $\Bb^j$:
\begin{enumerate}[(i)]
\item \( \bd(x) = \bd(y) = j\), for some \(0\le j \le d\), and

\item \(\bo^j(x) \etail \bo^j(y)\) in $\Bb^j$.
\end{enumerate}
That is, if one writes \(x\sim_j y\) to specify the border dimension, and call $\rel_\Bb^j$ the corresponding subrelation, then $\rel_\Bb$ is the union of the $\rel_\Bb^j$.
As a consequence of minimality, for any $j<d$, the set of paths of border dimension $j$ is dense in $\Pi_\infty^d$ (Proposition~\ref{bratteli10.prop-dense}), and it has measure zero with respect to any translation invariant measure \cite{RS98}.
So the relation $\rel_\Bb^j$, for $j<d$, is defined on a {\em thin set}.
The relation $\rel_\Bb^d$ is the standard $AF$-relation on $\Bb^d$.
But for $j<d$, it is important to notice that  $\rel_\Bb^j$ is {\em not} an $AF$-relation (see Remark~\ref{bratteli10.rem-comparisonER2}).
We prove here the following (Theorem~\ref{bratteli10.thm-eqrel} and Corollary~\ref{bratteli10.cor-etale}).

\vspace{.2cm}

\noindent {\bf Theorem} {\em The equivalence relation $\rel_\Bb$ is {\em \'etale}, and it is homeomorphic to $\rel_\Xi$.
}

\vspace{.2cm}

This gives a decomposition of $\Rr_\Xi$ in sub-equivalence relations:
the $AF$-equivalence relation homeomorphic to $\rel_\Bb^d$, and the ``missing parts'' which are pairs of tilings of border dimension smaller than $d$.

The definition of the topology of $\rel_\Bb$ is technical, and requires a finer analysis of the combinatorics of the substitution, as well as its encoding in the multi-diagram.
By minimality, a path $x$ can be the limit of a sequence \((x_n)_{n\in\NM}\) in $\Pi_\infty^d$, with \(\bd(x_n)\neq \bd(x)\) for all $n$.
So we loosen up the notion of paths to allow for changes in border dimensions.
We introduce {\em generalized paths} in the multi-diagram, which are infinite paths that have tails in $\Bb^j$ for some $j$ but can start in $\Bb^i$ with $i<j$ (see Section~\ref{bratteli10.ssect-gpath}).
For this purpose we define {\em escaping edges} which go from one vertex $v$ in $\Bb^i$ at depth $n$ to a vertex $u$ in $\Bb^{j}, j>i$, at depth $n+1$, whenever the face corresponding to $v$ lies in the interior of the substitute of that associated with $u$ (see Figure~\ref{bratteli10.fig-escaping}, and section~\ref{bratteli10.ssect-escape}).
This determines the topology of $\rel_\Bb$ from the combinatorics of the multi-diagram.

\subsection{Perspectives}
\label{bratteli10.ssect-perspectives}
In this paper, we present how to describe in a combinatorial way the equivalence relation on a substitution tiling space.
Several questions arise from this work.
First, the question of a generalization to groupoids arising from more general tilings.
Then, the question whether Bratteli multi-diagrams provide a model for a certain class of equivalence relations.
Finally, the implications of this construction on the $C^*$-algebraic level.

As far as the generalization to more general tilings is concerned, it seems easy to make our construction work for substitution tilings with tiles with ``wild boundaries''.
Our definition of a face of a tile (Definition~\ref{bratteli10.def-face}) is indeed very combinatorial.
We could remove the assumptions that tiles are $CW$-complexes, to cover cases where tiles have fractal boundaries for instance.

What about equivalence arising from general tilings, or from $\ZM^d$ actions?
A generalization to simplicial tilings seems reasonable.
To any transversal of a minimal tiling space with finite local complexity, it is possible to associate a Bratteli diagram.
The construction of the diagram relies on the construction of \emph{refined tesselations} (expanding-flattening sequences in the sense of~\cite{BBG06}).
At step $n$, build a tiling space $\Omega_n$ whose prototiles can be tiled by tiles of $\Omega$.
However, the topology of these tiles becomes increasingly complicated: there is \emph{a priori} no good geometrical notion of a \emph{face} of such tiles.
In the present paper, a deliberate choice was made to give a definition of faces which is as combinatorial as possible. The only geometric
ingredient in our definition of a face (Definition~\ref{bratteli10.def-face}) is the notion of dimension.
In the case of a simplicial tiling, the dimension is a combinatorial quantity: the non-empty intersection of $k$ distinct $d$-dimensional simplices defines an object of dimension $d-k+1$.
Simplices and Delone triangulation are used in particular for this reason in the work of Giordano, Matui, Putnam and Skau~\cite{GMPS10}.
If one relaxes what it means for two tiles to intersect, it seems reasonable to define faces and induced substitutions using their formalism of {\em well-separated tesselations}.
It could be expected to extend the results of the present article to general tilings, and in particular to groupoids arising from
minimal $\ZM^d$ actions on a Cantor set.

On the $C^*$-algebraic level, Bratteli diagrams were originally used to classify $AF$-algebras.
It would be interesting to see whether our multi-diagram and the structure of the equivalence relation $\rel_\Bb$ that it encodes could shed some light on the structure of the tilings $C^*$-algebras.

\paragraph{Acknowledgements}
We would like to thank Ian Putnam and Thierry Giordano for useful discussions.
The question whether the tiling's equivalence relation could be encoded in a purely combinatorial way was raised by I. Putnam.

\section{Substitution tilings}
\label{bratteli10.sect-stilings}

In this section, we briefly define the notions of tile, tiling, tiling space,
and canonical transversal arising from a substitution rule.
Given a tiling, there is a natural action on the associated tiling space, which
is given by translation.
The orbit equivalence relation induced by the action restricts to the
transversal.
We give some details on the topology of these equivalence relations.
Finally, we define \emph{decorations} of tiles and faces of tiles.


\subsection{Some notions on substitution tilings}
\label{bratteli10.ssect-stilings}

We work in the $d$-dimensional Euclidean space $\RM^d$.
Let us first define some vocabulary.
We refer the reader to \cite{Rob04,Frank08} for a complete exposition.

\begin{defini}
\label{bratteli10.def-tiling}
\begin{itemize}
 \item A {\em tile} is a compact subset which is homeomorphic to a ball.
 \item A {\em partial tiling} is a set of tiles $p = \{ t_i\}_{i\in I}$ which
have pairwise disjoint interiors.
We set \( \cup p := \bigcup_{i\in I} t_i\).
 \item A {\em patch} is a finite partial tiling;
 \item A {\em tiling} is a partial tiling with support $\RM^d$.
\end{itemize}
\end{defini}

We add the following hypothesis on tiles.

\begin{hypo}
\label{bratteli10.hyp-CW}
The tiles and tilings have a \emph{cellular} structure, that is:
\begin{enumerate}[(i)]
 \item tiles are assumed to be finite CW-complexes;
 \item in a (partial) tiling or a patch, the intersection of any number of tiles
is either empty or a sub-complex of each of them.
\end{enumerate}
\end{hypo}
This allows to define faces, and in particular to define faces of a certain
dimension. See Definition~\ref{bratteli10.def-face} in
Section~\ref{bratteli10.ssect-indsub}.

There is a natural action of $\RM^d$ on the set of tiles by translation.
This action extends to patches, partial tilings and tilings:
\[
 T + x := \{ t + x \ ; \ t \in T \} \quad \text{for } x \in \RM^d.
\]

Notice that tilings are not regarded up to translation: if $T$ is a tiling,
then $T$ and $T+x$ are different for $x  \neq 0$.
We do not consider tiles and patches up to translation either, but as subsets
of $\RM^d$.

\begin{defini}
A puncturing of the tiles is a function $\punc$, which associates to a tile $t$
a point in its interior, such that:
\[
 \forall x \in \RM^d \quad \punc (t + x) = \punc (t) + x.
\]
\end{defini}

Then, the tile $t$ is said to be \emph{punctured}, and $\punc (t)$ is called
the \emph{puncture} of $t$.
A set of punctured tiles is a set of tiles with a puncturing function defined
on it.

\begin{notat}
We write \(\spt(p)\) for the support of a patch $p$.
Given a tiling $T$, we set
\begin{equation}\label{bratteli10.eq-punc}
  T^\punc = \bigl\{ \punc(t) \ : \ t\in T \bigr\} \,,
\end{equation}
the set of punctures of its tiles.
If $t$ is a tile and $p$ a patch, the notations 
\[
t \in T \,, \quad \text{ and } \quad p \subset T\,,
\]
respectively mean ``$t$ is a tile and $p$ is a patch of the tiling $T$, at the
positions they have as subsets of $\RM^d$''.
We will use the following notation
\[
t \textrm{ appears in } p \,, \quad  t \textrm{ appears in }  T\,, \quad
\textrm{ and } \quad p \textrm{ appears in }  T \,,
\]
if there exists $a \in \RM^d$ such that we respectively have
\(t+a \in p\), \(t+a \in T\), and \(p+a \subset T\).
\end{notat}

Let us now define substitution tilings. Start with a set of prototiles, then
define a substitution rule on it.
\begin{defini}
\label{bratteli10.def-prototile}
\begin{itemize}
\item A \emph{prototile set} $\proto$ is a finite family of equivalence
classes of tiles of $\RM^d$ under translation.
\item A set of \emph{punctured prototiles} is a set of prototiles $\proto$,
together with a puncturing function $\punc$ defined on the set of all tiles
with class in $\proto$.
\end{itemize}
\end{defini}

By abuse of notation, we may identify an element $t \in \proto$ with its
unique representative $t_0$ which satisfies $\punc(t_0) = 0_{\RM^d}$.
We also say that a ``patch with tiles in $\proto$'' is a patch whose tiles have
their translational classes in $\proto$.
We define similarly a (partial) tiling with tiles in $\proto$.

Note that in the definition of tiles and prototiles, we allow ``labels'': it
is possible that two elements of $\Aa$ have the same shape, but a label
indicates that they should be regarded as different elements.


\begin{defini}\label{bratteli10.def-substitution}
 A substitution rule $\omega$ with inflation factor $\lambda$ on the
prototile set $\proto$ is a map which, to a tile $t$ of $\proto$, associates a
patch with tiles in $\proto$, such that:
\[
 \spt (\omega (t)) = \lambda \spt (t),
\]
and for all $x \in \RM^d$,
\[
 \omega (t+x) = \omega (t) + \lambda x.
\]
\end{defini}

\begin{defini}
Given a substitution $\omega$, define the \emph{Abelianization
matrix} of $\omega$ as the matrix $A = (A_{ij})_{i,j \in I}$, where
$I$ is in bijection with the set of prototiles \emph{via} a map $i \mapsto t_i$,
and such that
\[
 A_{ij} = \text{number of occurences of } t_i \text{ in } \omega(t_j).
\]
\end{defini}

The substitution allows to define what ``acceptable tilings'' (with respect to
$\omega$) are, and to define $\Omega$, the set of all acceptable tilings.

\begin{defini}
 The tiling space $\Omega$ associated to $\omega$ is the set of all tilings $T$
such that for all patch $p \subset T$, there exists $t \in \proto$ and $n \in
\NM$, such that $p$ appears in $\omega^n (t)$.
\end{defini}

It is clear that for any $T \in \Omega$ and $x \in \RM^d$, $T+x \in \Omega$.
Therefore, there is a natural $\RM^d$ action on $\Omega$.
It is classical that $\Omega$ is not empty: it is
possible to build a fixed point of some power of the substitution
$\omega$; such a fixed point then belongs to $\Omega$ (see
Figure~\ref{bratteli10.fig-subst} for an example).

We now make some assumptions on the substitution.
\begin{hypo}\label{bratteli10.hypo-primitive}
\begin{enumerate}[(i)]
  \item The substitution is \emph{primitive}: the associated
Abelianization matrix is primitive;
  \item the substitution is \emph{strongly aperiodic}: for all tiling
$T \in \Omega$, $\big( T+x = T \big) \Rightarrow \big( x = 0 \big)$;
  \item the tiling space has \emph{finite local complexity} (FLC): there
are finitely many patches of a given size, up to translation.
\end{enumerate}
\end{hypo}

Primitivity for a matrix $A$ means that there is some integer $n > 0$ such that
all entries of $A^n$ are strictly positive.

The tiling space $\Omega$ can be given a topology. It is defined by the
following basis.
For a patch $p$, and $r > 0$, let
\begin{equation}
\label{bratteli10.eq-openOmega}
 \Omega(p,r) :=  \{ T \in \Omega \ ; \ \exists x \in \RM^d, \nr{x} < r, \ 
 p \subset (T-x) \}.
\end{equation}
The sets $\Omega(p,r)$ form a basis for a topology on $\Omega$.
With this topology, the $\RM^d$-action by translations is continuous.

\begin{proposi}
 With Hypothesis~\ref{bratteli10.hypo-primitive}, $(\Omega,\RM^d)$ is a
compact and minimal dynamical system.
\end{proposi}
Minimality means that all orbits are dense. Minimality is actually equivalent to
the combinatorial condition of \emph{repetitivity} on the tiling: all patches
repeat ``often'' in a certain sense.
The critical condition to ensure minimality is the primitivity of the
substitution.

We now define a transversal $\Xi$ for the action of $\RM^d$ in $\Omega$.
\begin{defini}
\label{bratteli10.def-Xi}
The \emph{canonical transversal} of $\Omega$ (with respect to the $\RM^d$
action) is:
\[
  \Xi := \{ T \in \Omega \ ; \ 0_{\RM^d} \in T^\punc \}.
\]
\end{defini}


It is easily shown that the relative topology of $\Omega$ restricted to
$\Xi$ is given by the following basis of open sets:
given a patch $p$ such that $0_{\RM^d} \in p^\punc$, consider
\begin{equation}
\label{bratteli10.eq-openXi}
\Xi(p) = \{ T \in \Xi \ ; \ p \subset T \}\,.
\end{equation}

\begin{proposi}
The canonical transversal $\Xi$ is a Cantor set, that is a compact Hausdorff,
totally disconnected set, with no isolated points.
Furthermore, the sets defined in~\eqref{bratteli10.eq-openXi} form a basis
of clopen sets (sets which are both open and closed).
\end{proposi}

\subsection{Tiling equivalence relations and groupoids}
\label{bratteli10.ssect-gpd}

Let $\Omega$ be a tiling space, and let $\Xi$ be its canonical transversal.
We define two equivalence relations associated with $\Omega$ and $\Xi$.

\begin{defini}
\label{bratteli10.def-ERhull}
The equivalence relation of the tiling space is the set
\begin{equation}
\label{bratteli10.eq-ERhull}
\rel_\Omega = \bigl\{ \
(T,T') \in \Omega \times \Omega \ : \ \exists a \in \RM^d\,, \; T' = T + a \
\bigr\}
\end{equation}
with the following topology: a sequence $(T_n, T'_n = T_n + a_n)$ converges to $(T,T'=T+a)$ if \( T_n \rightarrow T\) in $\Omega$ and \(a_n \rightarrow  a\) in $\RM^d$.

The equivalence relation of the transversal is the restriction of $\rel_\Omega$ to $\Xi$: 
\begin{equation}
\label{bratteli10.eq-ERtrans}
\rel_\Xi = \bigl\{
(T,T') \in \Xi \times \Xi \ : \ \exists a \in \RM^d\,, \; T' = T + a \
\bigr\}
\end{equation}
\end{defini}

Note that the equivalence relations are {\em not} endowed with the relative
topology of \(\rel_\Omega \subset \Omega \times \Omega\) and \(\rel_\Xi \subset
\Xi \times \Xi\).
For example, by minimality, for $a$ large, $T$ and $T+a$ might be close in
$\Omega$, so that $(T,T+a)$ is close to $(T,T)$ for the relative topology, but
not for that from $\Omega \times \RM^d$.
The map \( (T,a) \mapsto (T, T+a)\) from  $\Omega \times \RM^d$ to $\Omega
\times \Omega$ is injective because $\Omega$ is strongly aperiodic (contains no
periodic points), and its image is $\rel_Omega$.
We have actually transfered the topology of $\Omega \times \RM^d$ to
$\rel_\Omega$ via this map.

\begin{defini}
\label{bratteli10.def-etale}
An equivalence relation $\rel$ on a compact metrizable space $X$ is called {\em \'etale} when the following holds.
\begin{enumerate}[(i)]

\item The set \(\rel^2 = \{((x,y),(y,z)) \in \rel \times \rel \}\) is closed in
$\rel \times \rel$ and the maps sending $((x,y),(y,z))$ in $\rel\times \rel$ to
$(x,z)$ in $\rel$,  and $(x,y)$ in $\rel$ to $(y,x)$ in $\rel$ are continuous.

\item The diagonal \( \Delta(\rel)=\{ (x,x) : x\in X\}\) is open in $\rel$.

\item The range and source maps \( r,s : \rel\rightarrow X\) given by \(r(x,y)=x, s(x,y)=y\), are open and are local homeomorphisms.

\end{enumerate}

A set $O\subset \rel$ is called an {\em $\rel$-set}, if $O$ is open, and $r\vert_O$ and $s\vert_O$ are homeomorphisms.
\end{defini}
The collection of $\rel$-sets forms a base of open sets for the topology of
$\rel$.
For this topology, it is proven in \cite{Kel95} that $\rel_\Omega$ and
$\rel_\Xi$ are \emph{\'etale} equivalence relations.

The tiling space $\Omega$ has a foliated space structure with leaves identified
to $\RM^d$ and Cantorian transversals \cite{BBG06}.
The holonomy groupoid $\Gamma_\Xi$ of the canonical transversal $\Xi$ to
$\Omega$ encodes essential dynamical and topological properties of $\Omega$.

A {\em groupoid} \cite{Re80} is a small category (the collections of objects and morphisms are sets) whose morphisms are all invertible.
A topological groupoid, is a groupoid $G$ whose sets of objects $G^0$ and morphisms $G$ are topological spaces, and such that the composition of morphisms \( G\times G \rightarrow G\), the inverse of morphisms \(G\rightarrow G\), and the source and range maps \(G\rightarrow G^0\) are all continuous maps.

Given an equivalence relation $R$ on a topological space $X$, there is a natural topological groupoid $G$ associated with $\rel$, with objects $G^0=X$, and morphisms $G=\{ (x,x') : x\sim_\rel x'\}$.
The topology of $G$ is then inherited from that of $\rel$.

\begin{defini}
\label{bratteli10.def-gpd}
The {\em groupoid of the tiling space} is the groupoid of $\rel_\Xi$, with set of objects $\Gamma_\Xi^0 = \Xi$ and morphisms
\begin{equation}
\label{bratteli10.eq-gpd}
\Gamma_\Xi = \bigl\{ (T,a) \in \Xi \times \RM^d \ : \ T+a \in \Xi \bigr\}\,.
\end{equation}
\end{defini}
There is also a notion of \emph{\'etale} groupoids \cite{Re80}.
Essentially, this means that the range and source maps are local homeomorphisms.
It can be shown that $\Gamma_\Xi$ is an \emph{\'etale} groupoid~\cite{Kel95}.

 \subsection{Faces and induced substitutions}
 \label{bratteli10.ssect-indsub}

In this section, we define faces of tiles, and describe how $\omega$ induces
substitutions on faces of any dimension.
A key point is that we need \emph{decorated} faces (or \emph{collared} faces).
The use of decorations is related to the notion of \emph{border forcing}
introduced by Kellendonk~\cite{Kel95}.
In their paper~\cite{AP98}, Anderson and Putnam used collared tiles to build
approximants of the tiling space, and describe the tiling space as an inverse
limit.
Bratteli diagrams can also be seen as an inverse limit construction which describes
the transversal. Therefore, decoration is also an essential feature.

Consider a set of prototiles $\proto$ of dimension $d$, and a substitution
$\omega$ on it. Let $\Omega$ be the associated tiling space.
Remember (Hypothesis~\ref{bratteli10.hyp-CW}) that tiles are CW-complexes.
Furthermore, the tilings in $\Omega$ are cellular in the sense that the
intersection of two adjacent tiles in a tiling is a subcomplex of both (this
is the analogue of meeting face-to-face for a tiling by polygons).

\begin{figure}[!h]  
\begin{center}
\includegraphics[width=9cm]{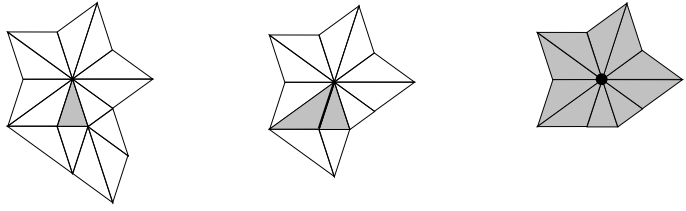}
\end{center}
\caption{\small{From left to right: faces of dimension $2$ ($p$ is a tile and $q$ its collar), dimension $1$ ($p$ is a pair of tiles), and dimension $0$ ($p=q$ in this case).}}
\label{bratteli10.fig-faces}
\end{figure}

\begin{defini}
\label{bratteli10.def-face}
Let $T \in \Omega$. A $j$-dimensional decorated face is
a pair of two patches of~$T$, $f := (p,q)$, satisfying the following
conditions:
\begin{enumerate}[(i)]
 \item $p$ and $q$ appear as subpatches in some tiling $T \in \Omega$;
 \item $\cap p := \bigcap_{t \in p}{t}$ is a $j$-dimensional (closed)
CW-complex;
 \item no intersection $(\bigcap_{t \in p}{t}) \cap t'$ with $t' \in T$
contains $\cap p$;
 \item $q$ is the set of all tiles of $T$ which intersect $\cap p$ (in
particular, $p$ is a subpatch of $q$).
\end{enumerate}
\end{defini}
See Figure~\ref{bratteli10.fig-faces} for an illustration in dimension $2$.

Given $f=(p,q)$ as above, the ``face'' itself is defined by the intersection of
all tiles of $p$, while its decoration is given by $q$.

\begin{notat}
Given $f =(p,q)$ a $j$-dimensional face, we set $\spt (f) := \cap p$ and
$\col (f) := q$.
\end{notat}

We extend the puncturing function to faces, such that $\punc (f) \in \spt(f)$,
in a coherent way with respect to translations: $\punc (f+x) = \punc(f) + x$.

\begin{defini}\label{bratteli10.def-protofaces}
For $0 \leq j \leq d$, define $\face{j}$ as the set of all equivalence classes
of $j$-dimensional (punctured) faces up to translation.
\end{defini}

All these sets are finite, by the finite local complexity property.
By abuse of notation, we may consider $f \in \face{j}$ as a specific
representant of such an equivalence class.

Note that as a particular case, a $d$-dimensional face is actually a
tile $t \in \proto$, together with a label (its collar).
This collar is the set of all tiles intersecting $t$ in a given tiling.
An example is given in Figure~\ref{bratteli10.fig-faces} (left).

We use the same terminology as for tiles: a $j$-dimensional face $f$
\emph{appears} in a tiling $T$ if there is some $a \in \RM^n$ such
that of $f + a$ is included in $T$ (in the sense that both patches defining
$f+a$ are included in $T$).
The face $f$ is in $T$ (noted abusively $f \in T$) if both patches defining $f$
are included in $T$ (at the same position).

The substitution $\omega$ extends naturally to all the sets $\face{j}$ as
follows.

\begin{defini}\label{bratteli10.def-indsubst}
Let $f = (p,q) \in \face{j}$.
Consider all pairs $(p_i,q_i)$, with $p_i$ a subpatch of $\omega (p)$ and
$q_i$ a subpatch of $\omega(q)$ which define a $j$-dimensional face in the sense
above.
The substitution $\omega_j (f)$ of $f$ is the set of all pairs $(p_i,q_i)$ which
satisfy $\cap p_i \subset \lambda (\cap p)$.
\end{defini}

\begin{figure}[htp]  
\begin{center}
\includegraphics[scale=0.6]{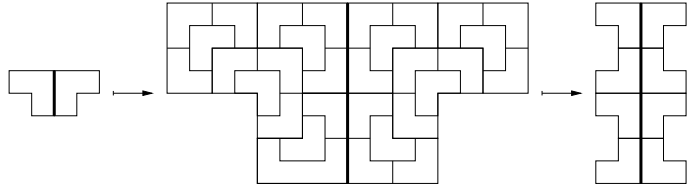}
\end{center}
\caption{\small{The substitution of the chair tiling induced on boundaries
(decorations not shown). On the left: the face defined by a pair of tiles; in
the middle: the substitution of the pair of tiles; on the right: the substitution
of the face (which is made of four pairs of tiles).}}
\label{bratteli10.fig-chair-bord-subst}
\end{figure}

Notice in particular that there is an induced substitution on decorated faces,
see Figure~\ref{bratteli10.fig-chair-bord-subst}.
The substitution on decorated tiles is primitive, and it can be proved that the
tiling space associated with it is conjugate to the tiling space $\Omega$.

Let give some precision on the definition of faces, to define intersections and
boundaries of faces.
\begin{defini}
 Let $f = (p,q)$ and $f' = (p',q')$ be two faces of arbitrary dimension.
We write that, $f$ and $f'$ \emph{intersect} if all the following conditions
hold:
\begin{enumerate}[(i)]
 \item $p \subset q'$;
 \item $p' \subset q$;
 \item $q \cup q'$ is a well-defined patch; that is the collars of $f$ and $f'$
match.
\end{enumerate}
%
\end{defini}

\begin{defini}\label{bratteli10.def-boundary}
A face $f' = (p',q')$ is \emph{on the boundary of} a face $f = (p,q)$ if the
following conditions both hold:
\begin{enumerate}[(i)]
 \item $p \subsetneq p'$ (which implies $\cap p' \subset \cap p$);
 \item $q' \subsetneq q$.
\end{enumerate}
\end{defini}

\begin{rem}
 In the above definition, if $f'=(p',q')$ is on the boundary of $f=(p,q)$, then
the dimension of $f'$ has to be strictly smaller than the dimension of $f$.
Indeed, by definition, for any $p''$ containing $p$, the dimension
of $\cap p''$ has to be strictly smaller than the dimension of $p$ (in the
sense of the dimension of a CW-complex).
\end{rem}

The fact that a (decorated) faces uniquely defines all its adjacent faces is
pictured in Figure~\ref{bratteli10.fig-hor}.
Notice also that the decoration of a face is included in the decoration of the
tile.

\begin{lemma}
Let $f = (p,q)$ and $f' = (p',q')$ be two distinct faces which intersect.
Then there is a face $f''$ of dimension strictly smaller than those of $f$ and
$f'$, such that $f''$ is on the boundary of both $f$ and $f'$.
\end{lemma}

\begin{proof}
First, remark that by the definition of the intersection of two faces, the
patch $p \cup p'$ is a well defined patch (it is actually a subpatch of both
$q$ and $q'$).
Let $p''$ be the subpatch of $q \cup q'$ such that $\dim(\cap p'') =
\dim(\cap (p \cup p'))$, and which is maximal for this property.
Note that $p''$ contains necessarily $p$ and $p'$.
Let $q''$ be the subpatch of $q \cup q'$ which contains all tiles intersecting
$\cap p''$. Since $\cap p''$ is included in both $\cap p$ and $\cap p'$, then
$q''$ is a subpatch of both $q$ and $q'$.
By the definitions, we built a face $f''$ which is on the boundary of both
$f$ and $f'$.
By the remark above, the dimension of this face is strictly smaller the
dimension of $f$ and $f'$, and the lemma is proved.
\end{proof}

\begin{notat}
We use the notation ``$f' \in \partial f$'' to say that $f'$ is on the boundary
of $f$, and $f'' \in f \cap f'$ to say that $f''$ is in the intersection of
$f$ and $f'$ in the above sense.
\end{notat}

The following technical lemmas will be used in the next sections.

\begin{lemma}
\label{bratteli10.lem-rho}
 There exists $\rho > 0$ such that:
\begin{enumerate}[(i)]
 \item for any face $f = (p,q) \in \face{j}$, a $\rho$-neighborhood of
$\cap p$ is included in the support of~$q$:
\begin{equation}\label{bratteli10.eq-rho}
 \big( \cap p + B(0,\rho) \big) \subset \big( \cup q \big).
\end{equation}
 \item If any two faces $f=(p,q)$ and $f'=(p',q')$, are $\rho$-close, then they
intersect:
\begin{equation}
 \big(d ( \cap p, \cap p') < \rho \big) \Rightarrow \big( f \cap f' \neq
\emptyset \big).
\end{equation}
\end{enumerate}
\end{lemma}

\begin{proof}
 Let $f \in \face{j}$. Assume $f$ is in some tiling $T \in \Omega$.
Assume that for all $\rho > 0$, there is some $x \in (\cap p + B(0,\rho))
\setminus (\cup q)$.
In particular, it is possible to chose a sequence $(x_n)_{n \in \NM}$ which
converges to a point $x \in \cap p$ but no $x_n$ belongs to $\cup q$.
By finite local complexity (and up to extraction of a subsequence), we may
assume that all $x_n$ belong to the same tile $t \in T$.
Then, since $t$ is closed, $\lim_n (x_n) \in (\cap p) \cap t$, and so by
definition of $q$, $t \in q$: it is a contradiction.
Therefore, there is a $\rho_f$ which satisfies \eqref{bratteli10.eq-rho} for
the tile $f$.
Take now $\rho = \min_{f}{\rho_f}$. It is positive by finite local complexity,
and satisfies \eqref{bratteli10.eq-rho} for any $f$. This proves (i).

Now, let $f, f'$ be two faces.
If $\cap p$ and $\cap p'$ are closer than $\rho$, then some point of $\cap p$ is
in the interior of $\cup q'$.
Therefore, all tiles of $p$ intersect the interior of $\cup q'$.
Therefore, $p \subset q'$.
Similarly, $p' \subset q$.
By definition, $f$ and $f'$ intersect.
\end{proof}

\section{Multi-diagram}
\label{bratteli10.sect-mdiag}

We define here a Bratteli multi-diagram associated with a substitution tiling of $\RM^d$.
We first recall the construction of the usual Bratteli diagram of substitution.

\subsection{Usual Bratteli diagram and $AF$-equivalence relations}
\label{bratteli10.ssect-diagAF}

Let us remind the reader how a Bratteli diagram is defined for a primitive
substitution~\cite{Kel95,For97}.

\paragraph{Construction of the diagram}

An example of a Bratteli diagram is given in the introduction (Figure~\ref{bratteli10.fig-BdiagFibo}).
The formal definition is the following.

\begin{defini}\label{bratteli10.def-brat-usual}
Let $\omega$ be a primitive and totally aperiodic substitution, with
prototile set $\proto$.
The stationary Bratteli diagram associated with $\omega$ is the graph
$\Bb = ( \Vv, \Ee)$, with
\begin{align*}
\Vv & = \left(\bigcup_{n \geq 1}{\Vv_n}\right) \cup \{\rac\}, &
\Ee & = \bigcup_{n \geq 0}{\Ee_n},
\end{align*}
where all the $\Vv_n$ are copies of $\proto^d$ (the set of decorated
prototiles, Definition~\ref{bratteli10.def-protofaces}), and there is an edge $e
\in \Ee_n$ ($n \geq 1$) between $v \in \Vv_n$ and $v' \in \Vv_{n+1}$ if and only
if there is an occurrence of the tile corresponding to $v$ in the substitution
(in its decorated version, see Definition~\ref{bratteli10.def-indsubst}) of the
tile corresponding to $v'$.
Finally, there is a single edge in $\Ee_0$ between the root $\rac$ and each
vertex of $\Vv_1$.
\end{defini}

The adjacency of edges and vertices is given by two maps $r$ and $s$ (range
and source maps), such that $r: \Ee_n \ra \Vv_{n+1}$ and $s: \Ee_n \ra \Vv_n$.

\begin{defini}
A \emph{path} in the Bratteli diagram $\Bb$ is a sequence of edges \(\gamma=
(e_n, \ldots, e_{m})\), for \(n<m\) and \(m\in \NM\cup\{\infty\}\), satisfying
$e_i \in \Ee_i$ and $r(e_i) = s(e_{i+1})$ for all $i $.
We denote by $\Pi_{n,m}$ the set of such paths.
If $m<\infty$ we extend the function $r$ to the $\Pi_{n,m}$, so that $r(\gamma) \in \Vv_m$
\end{defini}

We will use the shorthand notations $\Pi_n$ and $\Pi_\infty$ for $\Pi_{0,n}$ and $\Pi_{0,\infty}$ respectively
We endow each $\Ee_i$ with the discrete topology, and $\Pi_\infty$ with the
relative topology of the product topology on $\Pi_{n \geq 0}{\Ee_n}$.
Since the relation $r(e_i) = s(e_{i+1})$ is closed, it is clear that $\Pi_\infty$
is a compact and totally disconnected set (as a closed subset of a Cantor set).
The primitivity of the substitution ensures that it is itself a Cantor set.

\begin{notat}
For \(x \in \Pi_{n,m}\) and $n\le k<l \le m$ we denote by \(x_{[k,l]}\),
\(x_{[k,l)}\), \(x_{(k,l]}\), and \(x_{(k,l)}\) the restrictions of $x$ from
depths $k$ through $l$ with end points included or excluded.
For instance, if \(x\in \Pi_{\infty}\) we shall denote by \(x_{[n,\infty)}\)
the \emph{tail} of $x$ from depth $n$ on, and by \(x_{[0,n)}\) its head from the
root down to depth $n$ (excluded).
If \(\gamma, \eta\) are two paths with \(s(\gamma) = r(\gamma)\), we denote by
\(\gamma \cdot \eta\) the concatenated path
\end{notat}

A family of clopen sets which generates the topology can be given explicitly.
\begin{notat}
Given $\gamma \in \Pi_n$, with $n < \infty$, define:
\[
\vois{\gamma} := \{x \in \Pi_\infty \ ; \    x_{[0,n]} = \gamma \}.
\]
\end{notat}

\paragraph{$AF$-equivalence relation}

\begin{defini}
\label{bratteli10.def-AFER}
Let $\Bb$ be a Bratteli diagram and let
\[
\rel_m = \bigl\{
(x, \gamma) \in \rinfpath \times \Pi_{m} \, : \, r( x_{[0,m]}) = r( \gamma) \bigr\}\,.
\]
with the product topology (discrete topology on $\Pi_{m}$).

The {\em $AF$-equivalence relation} is the direct limit of the $E_m$ given by
\[
\rel_{AF} = \varinjlim_m \rel_m =
\bigl\{
(x,y) \in \rinfpath \times \rinfpath \, : \, \exists m, x_{[m,\infty)} = y_{[m,\infty)} \, \bigr\}\,,
\]
with the direct limit topology.
For $(x, y) \in \rel_{AF}$ we write $x\eaf y$, or $x \etail y$, and say that the paths are {\em tail equivalent}.
\end{defini}

It is well known that $\rel_{AF}$ is an $AF$-equivalence relation, as the direct limit of the compact {\it \'etale} relations $\rel_m$, see \cite{Phi05}.

\paragraph{The Robinson map}

We now define the Robinson map, which relates Bratteli diagrams and tiling
spaces. Let $\face{d}$ be the set of collared tiles.

\begin{figure}[htp]
\begin{center}
\includegraphics[scale=0.8]{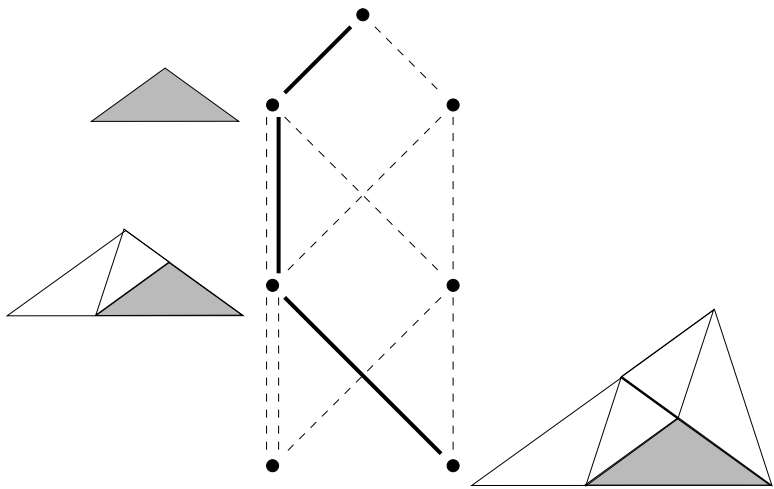}
\end{center}
\caption{\small{The Robinson map shown on a finite path for the Penrose tiling (decorations not shown).
For simplicity, the tiles are taken up to rotations and reflexions on this picture.}}
\label{bratteli10.fig-rob-map}
\end{figure}

For all $k$, let $\phi_0$ be the identification map $\Vv_k \ra \face{d}$.
The Robinson map $\phi$ is defined inductively, as a limit of maps $\phi_n$.
Fix $x = (x_0, x_1, x_2, \ldots) \in \Pi_{\infty}$ and proceed as follows.
\begin{itemize}
\item Define $\phi_1 (x)$ as the translate of the tile $\phi_0 (r(x_0))$ with
puncture at the origin;
\item Provided $\phi_{k-1}$ is defined, define $\phi_k (x)$ as the patch
$\omega^{k-1}(\phi_0(r(x_{k-1})))$, with the position of the origin determined
as follows:
since $e:=x_{k-1}$ encodes an occurrence of the tile $\phi_0 (s(e))$ in
$\omega (\phi_0 (r(e)))$, then it encodes an occurrence of the
patch $\omega^{k-2}(\phi_0(s(e)))$ in $\omega^{k-1} (\phi_0  (r(e)))$.
In other words, it encodes the inclusion of $\phi_{k-1} (x)$ in
$\phi_k (x)$. See Figure~\ref{bratteli10.fig-rob-map}
\end{itemize}

The properties of the $\phi_k$ make it possible to define $\phi (x)$ as the
union of the patches $\phi_k (x)$.
It can happen that this union is only a partial tiling.
However, since we used decorated tiles, there is a way to derive a full
tiling of $\RM^d$ from a decorated partial tiling.
We still call $\phi(x)$ this (undecorated) tiling.

\begin{defini}
\label{bratteli10.def-Robinpath}
Let $\phi : \Pi_\infty \lra \Xi$ be the map defined inductively from the
$\phi_n$.
This map is called the \emph{Robinson map}.
\end{defini}


\begin{theo}[\cite{Kel95}]
\label{bratteli10.thm-homeotrans}
The map $\phi : $ is a homeomorphism between the set of infinite rooted paths
in $\Bb$ and the canonical transversal of the tiling space.
\end{theo}

\subsection{Definition of the multi-diagram}
\label{bratteli10.ssect-mdiag}

We now turn to the definition of a generalized Bratteli diagram.
The basic idea is the following: some of the paths in a usual Bratteli diagram
only define partial tilings.
How much of these paths are there, and what is the structure of these special
tilings?
For one-dimensional tilings, it is known. There are finitely many such tilings,
and they correspond to fixed points of some power of the substitution.
In higher dimension, however, a \emph{nested structure} appears.
A half-tiling of the plane, for example, has a boundary which is a sequence
of one-dimensional faces.
In other words, it looks very much like a one-dimensional substitution tiling.
The generalized Bratteli diagrams contains the usual Bratteli diagram of the
$d$-dimensional substitution.
Also, for each dimension $0 \leq j \leq d-1$, it contains a Bratteli diagram
given by the substitution induced on $j$-faces.
It also has a horizontal structure: edges linking these different
diagrams. These edges encode the informations about faces being boundaries
of tiles.
Altogether, the information added to the original Bratteli diagram is purely
combinatorial and allows to define an equivalence relation on the set
of infinite paths.
This equivalence relation contains the tail-equivalence relation, and is mapped
homeomorphically \emph{via} the Robinson map to the translational equivalence
relation on the transversal of the tiling space, see
Section~\ref{bratteli10.sect-ERmdiag}.

\subsubsection{First step: dual diagram}

For technical reasons, it is more natural to start from a dual diagram rather
that from the usual Bratteli diagram.
The construction is done as follows.
Consider $\Bb_0$ the usual Bratteli diagram associated to $\omega$, as defined in Section~\ref{bratteli10.ssect-diagAF}.
Let us construct $\mbrat{d}=(\vtx^d, \edg^d)$ as follows.
\begin{itemize}
\item For all $n \geq 1$, the set $\vtx^d_n$ of vertices at depth $n$ in $\Bb^d$ is
isomorphic to the set of edges in $\Bb_0$ at depth $n$;
\item For all $n \geq 1$, there is an edge $e\in\edg^d_n$ between $s(e) \in \Vv_n^d$ and 
$r(e) \in \Vv_{n+1}^d$ if the corresponding edges are composable in $\Bb_0$;
\item Add a root and a set $\Ee_0^d$ of edges from the root to elements
of $\Vv_1^d$.
\end{itemize}
This new diagram is \emph{simple}: there is at most one edge between two given
vertices.
Therefore, a path in $\Bb^d$ (which is a sequence of composable edges)
is entirely given by the sequence of vertices it goes through.
Since vertices in $\Bb^d$ correspond to edges in $\Bb_0$, the map
\[
 (e_0, \ldots, e_n) \mapsto (r(e_0), \ldots, r(e_{n})).
\]
provides a canonical identification between paths in $\Bb_0$ and paths in
$\Bb^d$.
It makes therefore sense to define a Robinson map in this context: 
a finite path still corresponds to a partial tiling, and an infinite path
to a tiling.

\begin{figure}[htp]
\begin{center}
\psfrag{a}{{\small $(t_1,t'_1)$}}
\psfrag{b}{{\small $(t_2,t'_2)$}}
\psfrag{c}{{\small $(f,f')$}}
\psfrag{t}{{\small $t_1$}}
\psfrag{u}{{\small $\omega(t'_1)$}}
\psfrag{f}{{\small $f$}}
\psfrag{t'}{{\small $t_2$}}
\psfrag{u'}{{\small $\omega(t'_2)$}}
\psfrag{f'}{{\small $\omega(f')$}}
\psfrag{b1}{{\small $\Bb^2$}}
\psfrag{b0}{{\small $\Bb^1$}}
\includegraphics[width=12cm]{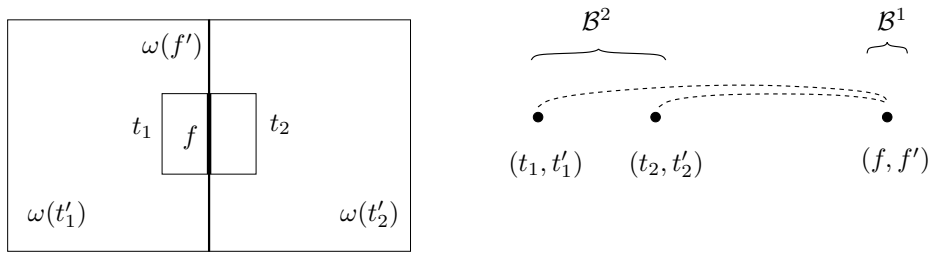}  
\end{center}
\caption{\small{Two vertices in the dual diagram $\Bb^2$ linked by horizontal edges (doted lines) to a common face in $\Bb^1$.}}
\label{bratteli10.fig-vertex-dual}
\end{figure}

What do vertices correspond to, \emph{via} the identification made by the
Robinson map?
Since a vertex of $\Bb^d$ correspond to an edge in $\Bb_0$, it correspond to
some tile, sitting inside a patch which is itself the substitution of a tile.
It is also possible to consider a vertex of $\Bb^d$ simply as a (decorated) tile,
with an additional label.
This additional label corresponds to the fact that this tile lies inside
of a given supertile in a predetermined position: it is not only information
about the neighborhood of the tile, but about its position in the hierarchical
structure of a tiling.
This is shown on the left of Figure~\ref{bratteli10.fig-vertex-dual}.

\subsubsection{Bratteli diagrams in lower dimensions and horizontal structure}
\label{bratteli10.ssect-hor}

For all $j < d$, let us build a Bratteli diagram associated to $\omega_j$.
We proceed exactly as before (Definition~\ref{bratteli10.def-brat-usual},
except that we do not include a root.
Then, we take the dual diagram.
The resulting diagram is called $\mbrat{j}$.

To sum up given $0 \leq j \leq d-1$, a vertex of $\mbrat{j}$, say
$v \in \Vv^{j}_n$, corresponds
to one occurrence of some $j$-face $f$ in the substitution of some other
$j$-face $f'$ under $\omega_j$.
There is an edge between $v$ and $v'$ if the inclusions are compatible
(that is $v'$ corresponds to an occurrence of $f'$ in some $\omega_j (f'')$).
This defines the diagram $\mbrat{j}=(\vtx^j,\edg^j)$.

\begin{figure}[htp]
\begin{center}
\includegraphics[width=.7\textwidth]{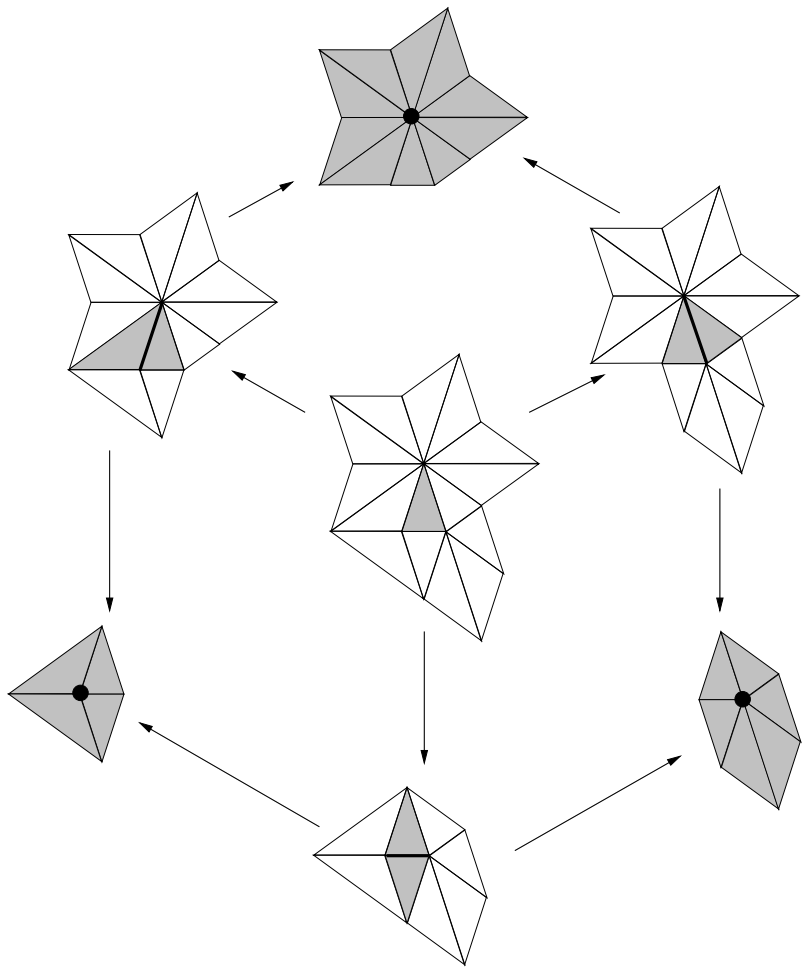}
\end{center}
\caption{\small{A tile (in grey in the center, with decoration shown around it) has here three $1$-dimensional faces (pairs of tiles),
and three vertices. The arrows shown here encode adjacencies.}}
\label{bratteli10.fig-hor}
\end{figure}

Each diagram is related to the other through the horizontal structure, which we
define now.
From the definition of decorated tiles and decorated $j$-faces, it is possible
to associate to any $j$-face its set of $(j-1)$-faces, independently of the
tiling in which they are sitting (see Figure~\ref{bratteli10.fig-hor}.
Indeed, the collar of a $(j-1)$-face is the set of all tiles which it
intersects.
Therefore, it is included in the set of all tiles which intersect
any $j$-face containing it, and the map which to a collared $j$-face associates
the set of $(j-1)$-faces on its boundary is well defined (remember
Definition~\ref{bratteli10.def-boundary}).

We shall denote by \(\Vv\) the union of the sets \(\Vv^j\) and \(\Ee\) the
union of the sets \(\Ee^j\) over $j=0, \ldots d$.
We now define horizontal edges, which links the diagrams $\Bb^j$ together.

\begin{defini}
For all $n \in \NM$, define the set of {\em horizontal edges} $\Hh_n$.
We have $\hor{n} = \bigcup_{j=1}^d{\hor[j]{n}}$, and there is an edge
$h \in \hor[j]{n}$ from $v \in \Vv^j$ to $v' \in \Vv^{j-1}$ if:
\begin{enumerate}[(i)]
 \item The vertex $v$ encodes the inclusion of some $g$ in $\omega_j (g')$;
 \item The vertex $v'$ encodes the inclusion of some $f$ in $\omega_{j-1}
(f')$,
with $f$ on the boundary of $g$ and $f'$ on the boundary of $g'$;
 \item The occurrence of $g$ in $\omega_j(g')$ encoded by $v$ actually lies on
the boundary and induces the inclusion encoded by $v'$.
\end{enumerate}
\end{defini}

We also set $\Hh$ to be the union of the sets $\hor[j]{}$ over $j=0, \dots d$.
We extend the range and source maps to $\hor[]{}$ as follows \( r : \hor[j]{n} \rightarrow \vtx^{j-1}_n\) and \( s: \hor[j]{n} \rightarrow \vtx^{j}_n\).

How horizontal structure corresponds to adjacency is represented on
Figures~\ref{bratteli10.fig-hor} and~\ref{bratteli10.fig-vertex-dual}.

\subsubsection{Escaping edges and vanishing faces}
\label{bratteli10.ssect-escape}

\begin{figure}[htp]
\begin{center}
\psfrag{b2}{{\small $\Bb^2$}}
\psfrag{b1}{{\small $\Bb^1$}}
\includegraphics[width=12cm]{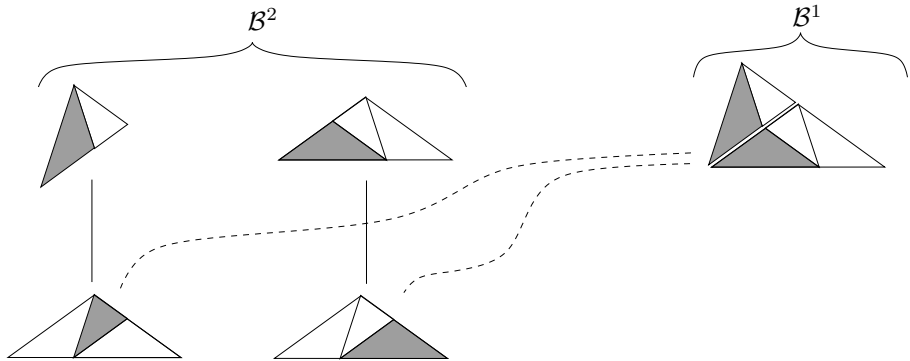}
\end{center}
\caption{\small{On the first level, the two tiles have a common face determined by their decorations (not shown here).
On the next level, this common face lies in the interior of a supertile, and is no longer visible on the diagram $\Bb^1$.
There is an escaping edge.}}
\label{bratteli10.fig-escaping}
\end{figure}

We complete the construction of the multi-diagram by constructing a new type of
edges, which we call \emph{escaping edges}.
The set of escaping edges keeps track combinatorially of the fact
that, in the substitution process, $\omega(t)$ contains faces in its interior,
which are ``created'' by $\omega$, and do not come from the substitution of
faces of $t$ (see Figure~\ref{bratteli10.fig-escaping}).
Any tile comes from the substitution of some other tile (which we call a
supertile in this definition in order to keep track of the hierarchy), but not
any face appears as a sub-face of a superface.

\begin{figure}[htp]
\begin{center}
 \psfrag{t1}{{\footnotesize $t_1$}}
 \psfrag{t2}{{\footnotesize $t_2$}}
 \psfrag{t3}{{\footnotesize $t_3$}}
 \psfrag{t4}{{\footnotesize $t_4$}}
 \psfrag{t5}{{\footnotesize $t_5$}}
 \psfrag{t6}{{\footnotesize $t_6$}}
 \psfrag{t7}{{\footnotesize $t_7$}}
 \psfrag{t8}{{\footnotesize $t_8$}}
 \psfrag{a}{{\footnotesize $f_2$}}
 \psfrag{b}{{\footnotesize $f_1$}}
 \psfrag{c}{{\footnotesize $f_4$}{\Large $\Bigg\lbrace$}}
 \psfrag{d}{{\Large $\Bigg\rbrace$}{\footnotesize $f_3$}}
 \psfrag{e}{{\footnotesize $f_5$}}
 \psfrag{1in3}[B1][B1][1][30]{{\footnotesize $t_1 \subset t_3$}}
 \psfrag{2in4}[B1][B1][1][30]{{\footnotesize $t_2 \subset t_4$}}
 \psfrag{3in5}[B1][B1][1][30]{{\footnotesize $t_3 \subset t_5$}}
 \psfrag{4in6}[B1][B1][1][30]{{\footnotesize $t_4 \subset t_6$}}
 \psfrag{5in7}[B1][B1][1][30]{{\footnotesize $t_5 \subset t_7$}}
 \psfrag{6in7}[B1][B1][1][30]{{\footnotesize $t_6 \subset t_7$}}
 \psfrag{7in8}{{\footnotesize $t_7 \subset t_8$}}
 \psfrag{ainc}[B1][B1][1][30]{{\footnotesize $f_2 \subset f_4$}}
 \psfrag{bind}[B1][B1][1][30]{{\footnotesize $f_1 \subset f_3$}}
 \psfrag{cine}[B1][B1][1][30]{{\footnotesize $f_4 \subset f_5$}}
 \psfrag{dine}[B1][B1][1][30]{{\footnotesize $f_3 \subset f_5$}}
\includegraphics[width=\textwidth]{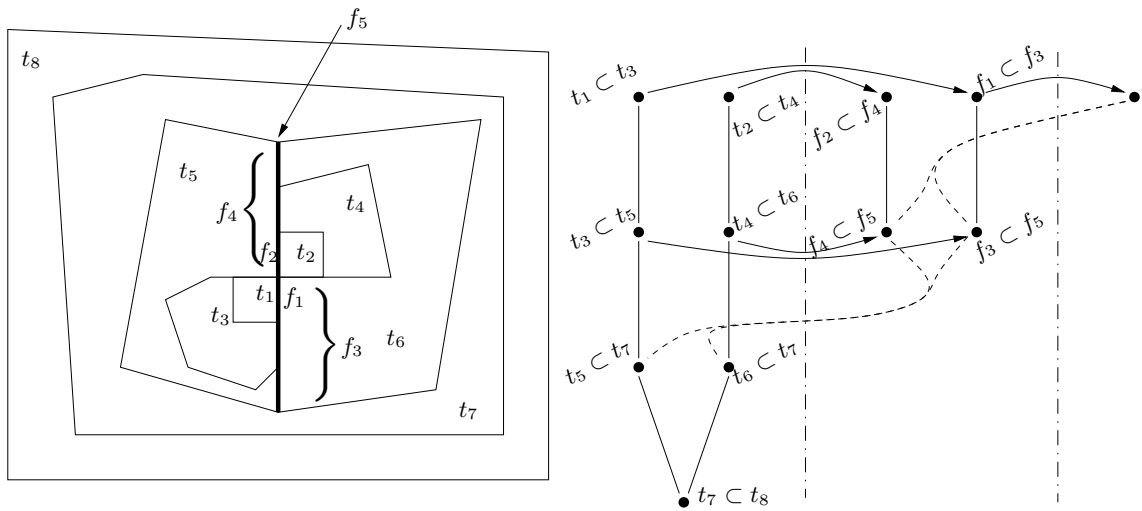}
\end{center}
\caption{\small{A situation in which two escaping edges follow one another.}}
\label{bratteli10.fig-escaping-double}
\end{figure}

\begin{defini}
Define the set of edges $\escaping$ on the multi-diagram as follows.
For all $n \geq 1$, there are edges $e \in \escaping$ in the following cases.
\begin{itemize}
\item From a vertex $v \in \Vv_n^{(j)}$ to a \emph{pair} of distinct vertices
$\{w,w'\} \subset \Vv_{n+1}^{(k)}$ in the case pictured in
Figure~\ref{bratteli10.fig-escaping}, that is if:
  \begin{itemize}
    \item $k > j$;
    \item $w$ and $w'$ correspond to inclusions of $k$-faces $g$ and $g'$ in the
\emph{same} $k$-superface $g''$;
    \item $v$ corresponds to the inclusion of some $j$-face $f$ in some
$j$-superface $f'$, and the intersection of $g$ and $g'$ contains $f'$.
  \end{itemize}
\item From a pair of vertices $\{v,v'\} \subset \Vv_n^{(j)}$ to a pair of
vertices $\{w,w'\} \subset \Vv_{n+1}^{(k)}$ in the case pictured in
Figure~\ref{bratteli10.fig-escaping-double}, that is if:
  \begin{itemize}
    \item $k > j$;
    \item $v$ and $v'$ correspond to inclusions of $j$-faces $f$ and $f'$ in the
\emph{same} $j$-superface $f''$;
    \item $w$ and $w'$ correspond to inclusions of $k$-faces $g$ and $g'$ in the
\emph{same} $k$-superface $g''$;
    \item $f''$ is the intersection of $g$ and $g'$.
  \end{itemize}
\end{itemize}
We extend the source and range maps to $\escaping$ as follows: for $e\in
\escaping$, $s(e)$ is a vertex or a pair of vertices in some $\vtx_n^j$, and
$r(e)$ is a vertex of a pair of vertices in some \(\vtx^k_{n+1}\), $k> j$.
\end{defini}

We finally define the Bratteli multi-diagram, by putting altogether the sets of
edges, horizontal edges and escaping edges.
We see this diagram as the union of the diagrams $\Bb^j$ over $j=0, \ldots d$,
linked together by the sets $\hor[]{}$ and $\escaping$ of horizontal and
escaping edges.
\begin{defini}
\label{bratteli10.def-multidiag}
The Bratteli multi-diagram associated with the substitution $\omega_d$ on $\Aa^d$ and its induced substitutions $\omega_j$ on $\Aa_j$ is the diagram \(\Bb = (\vtx, \edg, \hor, \escaping)\) with range and source maps $r,s,$ as defined previously on each sets of edges.
\end{defini}

\subsection{Paths in the multi-diagram}
\label{bratteli10.ssect-gpath}

We define several sets of paths in the diagram.
First, the paths in $\mbrat{d}$ are very similar to the paths defined in a
usual Bratteli diagram. We will allow non-rooted paths for technical reasons
and adapt the notations to this end.

We similarly define paths in $\mbrat{j}$, in the Bratteli diagrams
of lower dimensions.
Just like two vertices can be ``on the boundary of one another'' (if they are
linked by an horizontal edge), it is possible to define what it means for a
path to be on the boundary of another.
We will see in Section~\ref{bratteli10.ssect-border} what is the geometric
interpretation of this relation at the tiling level.

Finally, we define paths ``jumping up in dimension'': generalized paths on
the multi-diagram which use the escaping edges we defined above.

Since all our diagrams are simple, we index paths by vertices (or
maybe by vertices and pairs of vertices in the case of generalized paths).

\begin{notat}
We use the following notations:
$\chemin{n,m}{j}$ is the set of all paths in $\mbrat{j}$ starting at
depth $n$ and ending at depth $m$ (with possibly $m = \infty$).

Any $x \in \chemin{n,m}{j}$ is of the form
\[
 x = (v_n, v_{n+1}, \ldots, v_m),
\]
such that for all $i \in \{n, \ldots, m\}$, $v_i \in \Vv^j_i$, and there is an
edge $e \in \Ee^j_i$, with $s(e) = v_i$ and $r(e) = v_{i+1}$.

If $j=d$, then it could be that $n = 0$, in which case it means that the path
starts from the root, and we write $\Pi^d_\infty := \Pi^d_{0,\infty}$ for short.

Finally, $\Pi^j_{\bullet,m}$ with $m \in \NM \cup \{+\infty\}$ denotes the
union of $\Pi^j_{n,m}$ for $n \leq m$.
\end{notat}

There is a natural concatenation on paths. 
Given two paths $\gamma$ and $\eta$ with $r(\gamma)=s(\eta)$ we denote by $\gamma \cdot \eta$ their
concatenation.

Remember that a horizontal edge occurs when a tile is included in a supertile,
and this inclusion induces an inclusion of faces.

\begin{defini}
Let $\Ff$ be a set of vertices in some $\Vv^{j}$.
We define $\Ff'$ as the set of all vertices $v \in \Vv^{j-1}$ such that there
is an edge in $\Hh$ from an element of $\Ff$ to $v$:
\[
\Ff ' = r \left( \big( \restr{s}{\Hh} \big)^{-1}  (\Ff) \right).
\]
For all $0 \leq p \leq j$, we define $\Ff^{(p)} := (\Ff^{(p-1)})'$.
\end{defini}

In particular, if $x \in \chemin{\infty}{d}$ is a path, it is possible to
define all paths which lie in its boundary.
We adopt the following notation: if $x$ is a path, $\{x\}$ is the set of all
vertices it goes through, and so it makes sense to define $\{x\}'$.

\begin{defini}
\label{bratteli10.def-borderpath}
Given $x \in \chemin{\infty}{d}$ and $j \in \{0, \ldots, d\}$,
define:
\[
 \bo^k (x) = \{ y = (y_n, y_{n+1}, \ldots) \in \Pi^j_{n,\infty} \ ; \ n \in \NM \ 
; \ \forall i \geq n, \ y_i \in \{x\}^{(d-j)} \}.
\]
We define $\bo(x) = \bigcup_{0 \leq j \leq d}{\bo^j(x)}$.
We call $\bo^j (x)$ the set of \emph{boundary paths} of $x$ of dimension~$j$.
\end{defini}

In particular, $\bo^j (x)$ is a set of paths in $\Bb^j$. It could very well
be empty, and actually, we will see that it is generically empty for $j < d$
(see Remark~\ref{bratteli10.rem-comparisonER}).


\begin{defini}
The \emph{border dimension} of a path $x \in \Pi_{\infty}^d$ is the minimum
of all $j$ such that $\bo^j(x)$ is not empty 
\[
 \bd (x) = \min \{ j \le d \ ; \ \bo^j(x) \neq \emptyset \}
\]
\end{defini}

%
%

We now turn to the definition of generalized paths using escaping edges.

\begin{defini}
A \emph{generalized path} $x \in \gpath_{n,\infty}$ on the Bratteli multi-%
diagram is a sequence $(x_n, x_{n+1}, \ldots)$, where each $x_n$ is either a
vertex or a pair of vertices of $\Vv$, such that for all $i \geq n$, one of
these situations occurs:
\begin{itemize}
 \item there is an edge $e \in \Ee \cup \escaping$ such that $x_i = s(e)$ and
$x_{i+1} = r(e)$;
 \item if $x_i = (v,v')$ and $x_{i+1} = v''$, there is a pair of edges
$(e,e') \in (\Ee_i)^2$, such that $s(e) = v$, $s(e') = v'$ and $r(e) = r(e') =
v'';$.
\end{itemize}
\end{defini}


Let $x = (x_1, x_2, \ldots) \in \Pi^d_\infty$ (so that each $x_i \in \Vv_i$).
Let
\[
 \Ff_x = \{ v \in \Vv \ ; \  \exists  j \le  d\, \  v \in \{x\}^{(j)} \}.
\]
In other words, $\Ff_x$ is the set of all vertices of the path $x$, or in any
vertex set derived from $x$.

\begin{defini}
Let $x \in \Pi^d_{\infty}$.
A path $y =(y_n, y_{n+1}, \ldots ) \in \gpath_{n,\infty}$ is said to be \emph{derived from $x$} if
for all $k \geq n$
\begin{itemize}
 \item either $y_n$ is a single vertex and belongs to
$\Ff_x$, 

\item or $y_n$ is a pair of vertices, at least one of which belongs to
$\Ff_x$.
\end{itemize}
The set of all generalized paths derived for $x$ is:
\[
 \gbo (x) = \{ y \in \gpath_{n,\infty} \ ; \ n \in \NM \ 
\text{and y is derived from } x\}.
\]
\end{defini}

\subsubsection{Robinson map for generalized paths}

The Robinson map is generalized as follows, first for paths of any dimension,
then for generalized paths.

\paragraph{On paths}

Let $x = (x_m, x_{m+1}, \ldots) \in \Pi^{j}_{m, \infty}$.
We fix the notations for the proof: remember that a vertex $v$ in the
diagram corresponds to a face $f$ included in the substitution of a face
$f'$ (see Figure~\ref{bratteli10.fig-vertex-dual}).
By convention we say that $f$ is the face associated to the vertex $v$.
For all $n \geq m$, define $f_n = (p_n,q_n)$ the face associated to the vertex
$x_n$ (with $p_n$ possibly a single tile).
Remember that faces are given with a puncture, a distinguished point in their
support.

In the sequel, we will say that we ``put the origin in a patch at a certain
point'' to say that we take a translate of the patch so that this certain point
lies at the origin.
We also write $\omega(f)$ for the substitution of a face whose dimension is not
specified.

Define inductively $\phi_n$ on paths (with $n \geq m$) as:
\begin{enumerate}[(i)]
 \item $\phi_n (x)$ is a translate of $\omega^n (p_n)$;
 \item The position of the origin in $\phi_n (x)$ depends on
${x}_{[m,n]}$;
 \item $\phi_n (x)$ is a sub-patch of $\phi_{n+1}(x)$.
\end{enumerate}

Point (i) defines $\phi_n$. We just need to describe how to fix the origin in a
way which is compatible with point (iii).
The position of the origin in $\phi_n(x)$ is determined inductively as follows:
\begin{itemize}
 \item If $n=m$, put the origin in $\omega^n(f_n)$ at $\lambda^n \punc(f_n)$;
 \item Otherwise, since $p_{n-1} \subset \omega (p_n)$, the origin of $\phi_n
(x)$ is in the support of $\omega^{n-1} (p_{n-1})$, at a position determined by
$x_{[m,n-1]}$.
\end{itemize}

Condition (iii) above allows to define $\phi (x)$ as the union of the
$\phi_n (x)$.
Note that it could be a partial tiling.
However, using decorations, we show in~%
Proposition~\ref{bratteli10.prop-Robingpath} that it canonically extends
to a full tiling of $\RM^d$.

The collared version of $\phi$, noted $\phi^c$ is defined in the exact same way,
replacing the $p_n$ by the $q_n$.

\paragraph{On generalized paths}

What needs to be done in order to extend $\phi$ to any generalized
path $z \in \gpath_{m,\infty}$ is to describe what happens when some of the
$z_n$ are pairs of vertices.
If $z_n=(v_n,v'_n)$ is a pair of vertices, then $v_n$ corresponds to the
inclusion of a face $f'_n=(p_n,q_n)$ in the substitution of a face
$f''_n = (p''_n,q''_n)$, and $v'_n$ to the inclusion of $f_n = (p'_n,q'_n)$ in
the substitution of the {\em same} face $f''_n$.


Properties (i), (ii) and (iii) above still hold for $\grob_n$
whenever $z_n$ is a single vertex.
If $z_n$ is a pair of vertices, properties (ii) and (iii) are unchanged, and
property (i) becomes:
\begin{itemize}
 \item[($\tilde{\mathrm{i}}$)] $\phi_n (z)$ is a translate of $\omega_n
(p''_n)$.
\end{itemize}

The position is determined inductively by $z_{[m,n]}$ as follows.
\begin{itemize}
 \item If $n=m$ and $z_n$ is a single vertex the origin is determined as
above, 

\item if $n=m$ and  $z_n=(f_n,f'_n)$ is a pair of vertices put the origin in \(\omega^n(f''_n)\) at \(\lambda^n \punc(f''_n)\).

 \item If $z_n$ is a single vertex, then the face (resp. the faces) defined by
$z_{n-1}$ is (are both) included in $f_n$.
The origin is in the support of the faces defined by $z_{n-1}$, at a position
determined by $z_{[m,n-1]}$;
 \item If $z_n$ is a pair of vertices, the origin is in the support of $f_n \cap
f'_n$, at a position determined by $z_{[m,n-1]}$ (remember that $f_n \cap f'_n$
contains any face defined by the vertex or the vertices of $z_{n-1}$).
\end{itemize}

As for $\phi$, define $\grob^c_n$ like $\grob_n$, replacing the $p_n$ by the $q_n$.


\begin{proposi} 
\label{bratteli10.prop-Robingpath}
There is a continuous mapping, called Robinson map,
\begin{equation}
\label{bratteli10.eq-Robingpath}
\grob : \infgpath \lra \Omega\, ,
\end{equation}
such that for all $n \in \NM$, $\grob_n (z) \subset \grob (z)$.
\end{proposi}

\begin{proof}
By Lemma~\ref{bratteli10.lem-rho}, for all $n$, $q_n$ contains a
$\rho$-neighborhood of $p_n$.
Therefore, $\grob^c_n(x)$ contains a $(\rho \lambda^n)$-neighborhood of
$\grob_n(x)$.
In particular, since $\grob_n (x)$ contains the origin, $\grob^c_n (x)$
contains a $(\rho \lambda^n)$-neighborhood of $0_{\RM^d}$.
Therefore $\grob^c (x)$ is a tiling of all $\RM^d$.

\end{proof}

Note that $\grob$ is \emph{a priori} not one-to-one, and never onto.
Indeed if \(z\in\tilde{\Pi}_{m,\infty}\)
and \(z' = z_{[m',\infty)}\) for some $m'>m$, we might have
\(\grob^c_n(z) = \grob_n^c(z')\) for all $n\ge m'$ (if the punctures agree) and
therefore \(\grob(z)=\grob(z')\).

\begin{notat}
We use the following convention
\begin{equation}
\label{bratteli10.eq-notation1}
T_x := \rob(x)\,, \textrm{ for } x \in \rinfpath^d \,, \qquad \textrm{and} \qquad
x_T := \rob^{-1} (T)\,, \textrm{ for } T\in\Xi\,,
\end{equation}
where $\rob$ is the homeomorphism of Theorem~\ref{bratteli10.thm-homeotrans}.
And similarly we will set
\begin{equation}
\label{bratteli10.eq-notation2}
\tilde{T}_z := \grob(z)\,, \textrm{ for } z \in \infgpath^d \,.
\end{equation}
where $\grob$ is the map of Proposition~\ref{bratteli10.prop-Robingpath}.
\end{notat}

\subsection{Borders of a path and of a tiling}
\label{bratteli10.ssect-border}

We show here the geometrical meaning of Bratteli multi-diagrams for tilings.

We remind the reader that given \( x\in \rinfpath^d\), the set of infinite paths of $\bo(x)$ contained in $\Bb^i$ is called the {\em $i$--th border} of $x$, and denoted $\bo^i (x)$:
\[
\bo^i (x) = \bo(x) \cap \infpath^i \,.
\]
And $x$ is said to have {\em border dimension} $j$, and we write $\bd(x)=j$, if $\bo(x)$ contains an infinite path in $\Bb^j$, and no infinite path in $\Bb^i$ for any $i<j$:
\[
\bd(x)=j \quad \iff \quad \bo^j (x) \neq \emptyset \,, \ \textrm{and for all $j<i$\begin{large}\begin{Large}                                                                                               \end{Large}                                                                                  \end{large} }  \ \bo^i (x) = \emptyset \,.
\]


\begin{rem}
\label{bratteli10.rem-bdgpath}{\em
Any generalized path has a tail in $\infpath^j$ for some \(0\le j\le d\).
Hence \(x\in \rinfpath^d$ has border dimension $j$ if and only if there exists a generalized path in $\gbo(x)$ with tail in $\infpath^j$ and none with tail in $\infpath^i$ for any $i<j$.
}
\end{rem}

Note that,  if not empty, the set \(\bo^j (x)\) contains infinitely many paths in $\infpath^j$.
We show now that any two paths in \(\bo^j (x)\) are tail-equivalent, so that the class of \(\bo^j (x)\) is well defined.

\begin{lemma}
\label{bratteli10.lem-bdim}
Let \( x\in \rinfpath^d\) with $\bd(x)=j$, then any two infinite paths in \(\bo^j (x)\) are tail-equivalent in $\Bb^{j}$.
\end{lemma}
\begin{proof}
Let \(z,z' \in \bo^j(x)\).
Let $t_n$ be the tile corresponding to the vertex in $x$ at depth $n$, and \(f_n, f'_n \in \partial t_n\) the $j$-faces corresponding to the vertices in $z,z'$, respectively (at depth $n$ larger than some $m$ large enough so that they both are defined).
For all $n\ge n_0$ we have the inclusions 
\[
\lambda^{-n+m} \spt(f_{m}) \subset \spt(f_n) \subset \partial t_{n}\,, \quad \textrm{and} \quad 
\lambda^{-n+m} \spt(f'_{m}) \subset \spt(f'_n) \subset \partial t_{n}\,.
\]
The faces $f_m$ and $f'_m$ belongs to $\partial t_m$, and \(\lambda^{-n+m} t_m\) shrinks to a point as $n$ tends to infinity, 
therefore \(\dist(f_n, f'_n) \rightarrow 0\) as \(n\rightarrow \infty\).
Hence by Lemma~\ref{bratteli10.lem-rho} the faces must eventually intersect:
 \(f_n \cap f'_n\ne\emptyset\) for all $n$ greater than some $n_1$.

We now show that \(f_{n}=f'_{n}\) for all $n$ greater than some $n_2\ge n_1$, which proves \(z\etail z'\). 
It suffices to show that \(f_{n_2}=f'_{n_2}\).
Indeed if this holds, then both $f_{n_2+1}$ and $f'_{n_2+1}$ contain \(f_{n_2}=f'_{n_2}\) in their substitute, and so must be the same face of $t_{n_2+1}$.
And by immediate induction we get \(f_{n}=f'_{n}, n\ge n_2\).

Assume that it is not the case: for each  $n\ge n_1$ , \(f_n \cap f'_n\ne\emptyset\) but \(f_n \ne f'_n\).
Then there exists a $k_n$-face \(g_n \subset \partial f_n \cap \partial f'_n\), with $k_n<j$.
Since \(f_n \in \omega(f_{n+1})\) and  \(f'_n \in \omega(f'_{n+1})\) we have \(g_n \in \omega(g_{n+1})\) for all $n\ge n_1$.
If $g_{n_1}$ is a $k_{n_1}$-cell, then $g_n$ is a $k_n$-cell with $k_n\ge k_{n_1}$.
The sequence of cells dimensions $(k_n)_{n\ge n_1}$ is thus non decreasing and takes on values in the finite set \( S=\{k_{n_1}, k_{n_1}+1, \cdots j-1\}\).
Therefore it is eventually constant: there exists $k \in S$ and $n_2\ge n_1$, such that $k_n = k$ for all $n\ge n_2$.
In other words, the sequence of faces \((g_n)_{n\ge n_1}\) defines a generalized path in $\gbo(x)$ with tail in $\infpath^k$, for some $k\le j-1$.
Hence, by Remark \ref{bratteli10.rem-bdgpath}, we deduce that \(\bd(x) \le k \le j -1\), so \(\bd(x) \ne j\) and this is a contradiction.
\end{proof}

We can now characterize the notion of border dimension.
\begin{lemma}
\label{bratteli10.lem-limpatch}
Let $x \in \rinfpath^d$ be a rooted infinite path in $\Bb^d$, then
\[
\bd(x) = d \quad \iff \quad  \lim_{n\rightarrow + \infty} 
\dist\bigl( \rob_1(x), \partial \rob_n(x) \bigr) = + \infty \,,
\]
or equivalently
\[
\bd(x) < d \quad \iff \quad  \lim_{n\rightarrow + \infty} 
\dist\bigl( \rob_1(x), \partial \rob_n(x) \bigr) < + \infty\,,
\]
where \( \partial \rob_n(x)\) is a shorthand notation for the boundary of the support of $\rob_n(x)$.
\end{lemma}
\begin{proof}
Since the sequence of distances has always a limit in  \(\RM_+\cup\{+\infty\}\), the equivalence between the two statements is clear.
We prove the second statement.
Write \( d_n = \dist\bigl( \rob_1(x), \partial \rob_n(x) \bigr)\).
Let $t_n$ denote the tile corresponding to the vertex in $x$ at depth $n$.

If \(\bd(x)=j < d\), there exists an infinite path $z$ in some \(\bo(x) \cap \Pi^j_{m.\infty}\).
Let $f_n$ denote the $j$-cell corresponding to the vertex of $z$ at depth $n\ge m$.
For all $n\ge m$, \(f_n \in t_n\), so \(\lambda^n \spt(f_n)\) appears on the boundary of $\spt(\rob_n(x))$, thus
\(\dist\bigl( \rob_{m}(x), \partial \rob_n(x) \bigr) = 0\), and therefore
\(d_n = d_{m}\) for all $n\ge m$.

Conversely, by Lemma~\ref{bratteli10.lem-rho}, if $\rob_n(x)$ belongs to the interior of $\rob_{n+1}(x)$ then \(d_{n+1} > d_n + \rho\).
Otherwise if \(\partial \rob_n(x) \cap \partial \rob_{n+1}(x) \ne \emptyset\), then \(d_n=d_{n+1}\).
Therefore if the sequence $(d_n)_{n\in\NM}$ converges, it must be eventually constant.

Now assume that the sequence converges, say to $a\in \RM^d$, and let $n_1 \ge m$ be such that $d_n=a$ for all $n\ge n_1$.
We thus have \(\partial \rob_{n_1}(x) \cap \partial \rob_n(x) \ne \emptyset\) for $n\ge n_1$.
Let then $f'_n$ be a cell that appears in $t_n$ and such that \(\lambda^n f'_n \subset \partial \rob_n(x) \cap \partial \rob_{n+1}(x)\), $n\ge n_1$.
We therefore have \(f'_n \in \omega(f'_{n+1})\).
If $f'_{n_1}$ is a $k_{n_1}$-cell, then $f'_n$ is a $k_n$-cell with $k_n\ge k_{n_1}$.
The sequence of cells dimensions $(k_n)_{n\ge n_1}$ is non decreasing and takes on values in the finite set \( S=\{k_{n_1}, k_{n_1}+1, \cdots d-1\}\).
Therefore it is eventually constant: there exists \(k\in S\) and $n_2 \ge n_1$ such that $k_n = k$ for all $n\ge n_2$.
In other words, the sequence of faces \((f'_n)_{n\ge n_1}\) defines a generalized path in $\gbo(x)$ with tail in $\infpath^k$, for some $k\le d-1$.
Hence, by Remark \ref{bratteli10.rem-bdgpath}, we deduce that \(\bd(x) \le k \le d-1\), and so \(\bd(x)<d\).
\end{proof}

An easy consequence of minimality gives the following.
\begin{proposi} 
\label{bratteli10.prop-dense}
The set \( \bigl\{ x\in \rinfpath^d \, : \, \bd(x)=j \bigr\}\) is dense in $\rinfpath^d$, for any $j\le d$.
\end{proposi}
\begin{proof}
Call the above set $B_j$.
Pick \(x\in \rinfpath^d\), and \(y \in B_j\).
For $n\in \NM$, by minimality, there exists  $m\ge n$ such that the tile corresponding to the vertex in $x$ at depth $n$ appears in the $(m-n)$--th substitute of the tile corresponding to the vertex in $y$ at depth $m$.
This means that there exists a path \(\gamma \in \Pi_{n,m}^d\) with \(s(\gamma) = r(x_{[0,n]})\) and \(r(\gamma) = s(y_{[m,\infty)})\).
Define \( x_n = x_{[0,n]}  \cdot \gamma \cdot y_{[m,\infty)}\).
For all $n$, $x_n$ belongs to $B_j$, and the sequence $(x_n)_{n\in\NM}$ converges to $x$ in $\rinfpath^d$.
This proves that $B_j$ is dense in $\rinfpath^d$.
\end{proof}

\begin{figure}[h!]
 \begin{center}
   \includegraphics[width=10cm]{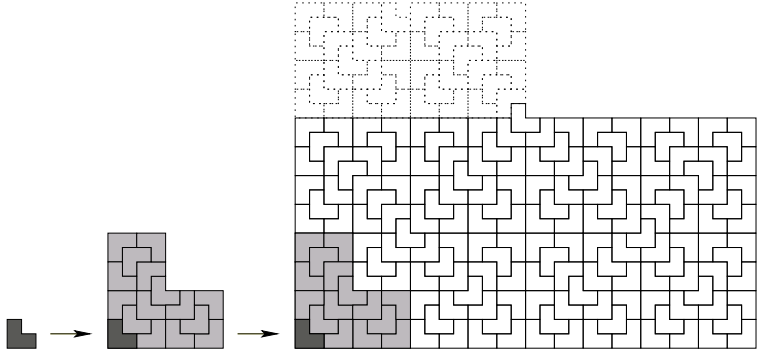}
   \caption{{\small A tiling of border dimension $0$ can be obtained as a fixed point of this map.}}
   \label{bratteli10.fig-bdim0}
 \end{center}
\end{figure}

\begin{figure}[h!]
 \begin{center}
   \includegraphics[width=10cm]{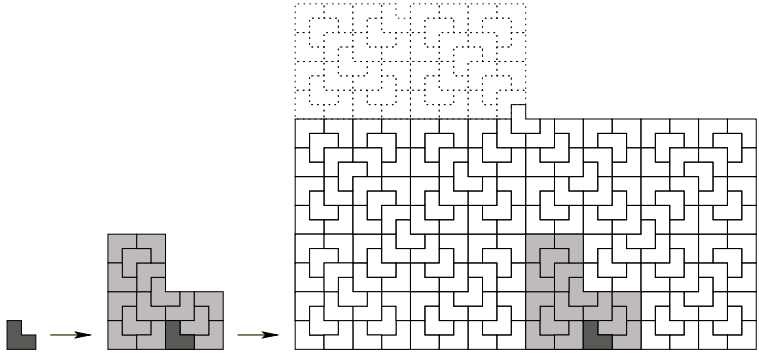}
   \caption{{\small A tiling of border dimension $1$ can be obtained as a fixed point of this map.}}
   \label{bratteli10.fig-bdim1}
 \end{center}
\end{figure}

\begin{figure}[!h]  
\begin{center}  
\includegraphics[width=5cm]{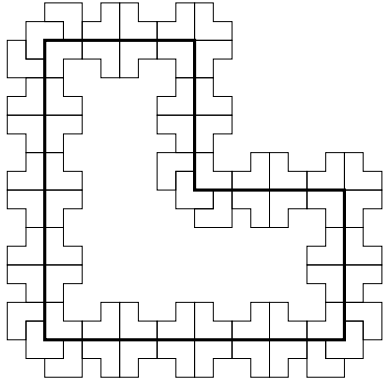}  
\end{center}
\caption{\small{$\lambda^2 \omega_2^{-2}(t)$ (thick line) and $ \omega_1^2( \omega^{-2}(t))$ for the chair tile $t$.}}
\label{bratteli10.fig-border}
\end{figure}

We show now that the border of a path corresponds a natural ``border'' of its associated tiling.
Figures~\ref{bratteli10.fig-bdim0} and~\ref{bratteli10.fig-bdim1} illustrate with the chair tiling how to obtain tilings of border dimension $0$ and $1$ respectively.
The following notion of border of a tiling was first introduced by Matui~\cite{Mat06} for substitution tilings of $\RM^2$.
\begin{defini}
\label{bratteli10.def-bordertiling}
For a tiling $T$ in $\Omega$, and \(0\le j\le d\), let \(T^j\) denote the union of the supports of the $j$-faces of $T$.
The {\em $j$-th border} of a tiling $T$ is defined as
\begin{equation}
\label{bratteli10.eq-bordertiling}
\bo^j(T) = \bigcap_{n\in\NM} \lambda^n \; \omega_d^{-n}(T)^j \,.
\end{equation}
\end{defini}
Remark that $\bo^j(T)$ is the ``support'' of the union of faces \(\bigcap_{n\in\NM} \omega_j^n \bigl( \omega_d^{-n}(T) \bigr)\), see Figure~\ref{bratteli10.fig-border}.

We have the elementary following properties.
\begin{lemma}
\label{bratteli10.lem-bordertiling}
\begin{enumerate}[(i)]

\item For any $a \in \RM^d$, one has \(\bo^j(T+a) = \bo^j(T) + a\). 

\item Let $T \in \Xi$.
For $n\in\NM$ let $t_n$ be the tile of \(\omega^{-n}(T)\) that contains the origin, and $p_n^j$ the union of the supports of the $j$-faces that intersect it.
One has
\[
\bo^j (T) = \bigcup_{m\in\NM}\bigcap_{n\ge m} \lambda^n \; p_n^j\,.
\]
\end{enumerate}
\end{lemma}
\begin{proof}
By definition of the substitution one has \(\omega_d^n(T+a) = \omega_d^n(T) + \lambda^n a\) and (i) follows immediately.

Let \(c^j(T)\) be the right hand side of the equation in (ii).
For any $T$ one has \(\lambda \omega_d^{-1} (T)^j \subset T^j\), so one can rewrite \(\bo^j(T)\) as \(\cup_{m\in\NM}\cap_{n\ge m} \lambda^n \omega^{-n}(T)^j\).
Since \(p_n^j \subset \omega_d^{-n}(T)^j\) one has \(c^j(T) \subset \bo^j(T)\).
To prove the other inclusion consider \(a \in \bo^j(T)\), so \(a \in \lambda^n \omega^{-n} (T)^j\) for all $n$.
Let $p_n$ denote the union of faces whose support is $p_n^j$.
Then $p_n$ contains $t_n$ in the interior of its support, so as \(n\rightarrow \infty\), $\lambda^n p_n$ eventually covers $\RM^d$.
Hence for all $n$ large enough \(a\in \lambda^n p_n^j\), and this proves the other inclusion \(\bo^j(T) \subset c^j(T)\) and completes the proof. 
\end{proof}

\begin{proposi} 
\label{bratteli10.prop-bordertiling2}
For any \(x\in \rinfpath^d\) one has
\[
\bd(x)=j \quad \iff \quad \bo^j (T_x) \ne \emptyset\,, \textrm{ and } \bo^i (T_x)= \emptyset\,, \ \forall i < j\,.
\]
\end{proposi}
\begin{proof}
One writes $t_n$ for the tile of \(\omega^{-n}(T)\) that contains the origin (and corresponds to the vertex of $x_T$ at depth $n$).

Consider the case $j<d$ first.
Assume \(\bd(x)=j\) and let \(z \in \bo^j(x)\). 
Let $f_n$ be the $j$-face on the boundary of the tile $t_n$, for $n$ larger than some $m$ so that it is defined.
For all $n\ge m$ we have \(\spt(f_n) \in p_n^j\), and by Lemma~\ref{bratteli10.lem-bordertiling} (ii) \(\bo^j (T) \ne \emptyset\).
Now if \(\bo^i (T) \ne \emptyset\) for some $i<j$, by Lemma~\ref{bratteli10.lem-bordertiling} (ii), for all $n$ larger than some $n_1$ there exist $i_n$-faces $f_n$ such that \(f_n \in \col(t_n)^i\) and with \(f_n\in\omega(f_{n+1})\), for some non decreasing integers $i_n$ taking values in the finite set \( S=\{i, i+1, \cdots j-1 \}\).
From the sequence of $i_n$-faces \((f_n)_{n\ge n_1}\) one builds a generalized path in \(\gbo(x)\) with tail in $\infpath^k$ for some $k\in S$.
By Remark \ref{bratteli10.rem-bdgpath} we deduce \(\bd(x)\le k \le j-1\), and hence \(\bd(x) \ne j\) which is a contradiction.

Conversely, if \(\bo^j (T_x) \ne \emptyset\) and \(\bo^i (T_x) = \emptyset\) for all $i<j$, then by Lemma~\ref{bratteli10.lem-bordertiling} (ii) there exists and infinite path in \(\bo(x) \cap \infpath^j\), and none in \(\bo(x) \cap \infpath^i\).
This proves that \(\bd(\bo)=j\).

Now consider the case $j=d$.
Assume \(\bd(x)=d\).
For all $n$, \(t_n \subset p_n^d = \col(t_n)^d\), and therefore \(\bo^d (T_x) \ne \emptyset\).
If for some $i<d$, \(\bo^i (T_x) \ne \emptyset\) then by the same argument given for the case $j<d$ above, one can build a generalized  path in \(\gbo(x)\) with tail in \(\infpath^i\) and this contradicts \(\bd(x)=d\).
For the converse, the proof is the same as for the case \(j<d\) above. 
\end{proof}

\begin{rem}
\label{bratteli10.rem-comparisonER}{\em
The $AF$-equivalence relation $\rel_{AF}$ on $\rinfpath^d$ induces an $AF$-equivalence relation $\rel'_{AF}$ on $\Xi$ via the Robinson map (Definition \ref{bratteli10.def-Robinpath} and Theorem~\ref{bratteli10.thm-homeotrans}): for $T\in \Xi$ we set
\[
[T]_{AF} = \bigl\{ T -a \in \Xi \, : \, x_{T-a} \eaf x_T \bigr\}\,.
\]
It is easy to see from Lemma~\ref{bratteli10.lem-limpatch} and Proposition~\ref{bratteli10.prop-bordertiling2}
\[
[T]_{AF} = [T]_{\rel_\Xi} \quad \iff \quad \bd(T) = d\,,
\]
and equivalently
\[
[T]_{AF} \subsetneq [T]_{\rel_\Xi} \quad \iff \quad \bd(T) < d \,.
\]
For example, the tilings shown in Figures~\ref{bratteli10.fig-bdim0} and ~\ref{bratteli10.fig-bdim1} give rise to only partial tilings of the plane (upper right quadrant, and upper half plane respectively), but extend uniquely to tilings of the whole plane {\it via} the Robinson map $\phi$, provided that the Bratteli diagram is built using decorated tiles.
However, the $\rel'_{AF}$ orbits of those tilings do not match their $\rel_\Xi$ orbits: the $AF$-relation cannot identify two tilings whose punctures lie on two different sides of the borders.
The authors explained this in detail in the first paper \cite{BJS10}, Section 3.3 (Remark 3.14 in particular), in a more general setting (without a substitution).

By proposition~\ref{bratteli10.prop-dense}, the set of tilings of border dimension less than or equal to $d-1$ is dense in $\Xi$.
A result of Radin and Sadun (``Property F'' in \cite{RS98}) shows however that this set has measure zero with respect to any invariant measure on $\Xi$.
Such a set is called {\em thin} in the literature.
So one sees that $\rel_\Xi$ differs from the $AF$-relation $\rel'_{AF}$ on thin set.
}
\end{rem}

\section{Equivalence relation in a Bratteli multi-diagram}
\label{bratteli10.sect-ERmdiag}

As a consequence of lemma \ref{bratteli10.lem-bdim}, we can associate to a path \(x\in \rinfpath^d\) of border dimension $\bd(x)=j$, the tail-equivalence class of its $j$--th border.

\begin{defini}
\label{bratteli10.def-eqrel}
We say  that two infinite paths $x$ and $y$ in $\rinfpath^d$ are {\em border equivalent}, and we write \(x \sim y\), if 
\begin{enumerate}[(i)]
\item \(\bd(x)= \bd(y)=j\) for some \(0\le j \le d\), and

\item \( \bo^j(x) \etail \bo^j(y)\) in $\Bb^j$.
\end{enumerate}
\end{defini}

This defines an equivalence relation on \(\rinfpath^d\): $\sim$ is clearly symmetric and reflexive, and transitivity follows from that of $\etail$.

\begin{rem}
\label{bratteli10.rem-eqrel}{\em
Any generalized path has a tail in \(\infpath^j\) for some $j$.
Hence two infinite paths $x$ and $y$ in $\Pi_\infty^d$ are {\em border equivalent} if and only if 
\[
\gbo(x) \cap \gbo(y) \cap \gpath_{\bullet,\infty} \neq \emptyset\,.
\]
}
\end{rem}

\begin{defini}
\label{bratteli10.def-ERlink}
Define the equivalence relation on $\rinfpath^d$
\[
\rel_\Bb = \bigl\{ (x,y) \in \rinfpath^d \times \rinfpath^d \ : \ x \sim y \bigr\} \,,
\]
with the following topology: \((x_n,y_n)_{n\in\NM}\) converges to \((x,y)\), if and only if
\begin{enumerate}[(i)]
\item \(x_n \rightarrow x\) and \(y_n \rightarrow y\) in $\rinfpath^d$, and

\item there exists \(m\in \NM\) such that \(\gbo(x_n) \cap \gbo(y_n) \cap \tilde{\Pi}_{m, \infty} \neq \emptyset\) for all $n$ large enough.

\end{enumerate}
\end{defini}


\begin{rem}
\label{bratteli10.rem-comparisonER2}{\em
We can write \(x \sim_j y\), in the case where \(\bd(x)=\bd(y)=j\) to specify the dimension of the borders, and call $\rel_\Bb^j$ the corresponding subrelation.
We thus view the equivalence relation $\rel_\Bb$ as the union of the $\rel_\Bb^j$ over \(j=0, \cdots d\):
\[
\rel_\Bb \; = \; \bigcup_{j=0}^d \; \rel_\Bb^j \,.
\]
By definition, if $\bd(x)=\bd(y)=d$ we see that $x\sim y$ if and only if $x \eaf y$, so that $\rel_\Bb^d = \rel_{AF}$.
So the equivalence relation $\rel_\Bb$ contains $\rel_{AF}$ as a natural subrelation.
However, by Remark~\ref{bratteli10.rem-comparisonER}, the inclusion is not an equality:
\[
\rel_{AF} = \rel_\Bb^d \subsetneq \rel_\Bb\,.
\]
Indeed the two relations only coincide on the set of paths or border dimension $d$.
They differ on the dense set of paths of border dimension less than $d-1$ (Proposition~\ref{bratteli10.prop-dense}).

Also, it is important to notice that $\rel_\Bb^j$ is {\em not} an $AF$-relation for $j<d$.
Indeed, its base (namely the set of paths of border dimension $j$) is neither compact, nor locally compact in $\rinfpath^d$.
}
\end{rem}

We now characterize $\rel_\Bb$.
We first state a technical lemma which proves that faces ``forces their borders'' (see~\cite{Kel95} for the original definition for tiles).
\begin{lemma}
\label{bratteli10.lem-facesborderforcing}
Let \(z \in \gpath_{m,\infty}\) and \(x \in \rinfpath^d\) such that \(z \in \gbo(x)\).
There exists an integer $k$ (that does not depend on $z,x$) such that for all $n\ge m$, \(\rob_n^c(x)\) appears in \(\grob_{n+k}^c(z)\). 
\end{lemma}
\begin{proof}
Let $t_n$ (respectively $f_n$) denote the tile (respectively face)  corresponding to the vertex in $x$ (respectively $z$) at depth $n\ge m$.
Since $z\in \gbo(x)$, for all $n\ge m$, $t_n$ appears in $\col(f_n)$.
Hence by Lemma~\ref{bratteli10.lem-rho}, \(\col(t_n)\) appears in \(\omega^{k}(\col(f_{n+k}))\) for any $k$ large enough, see Figure~\ref{bratteli10.fig-facesborderforcing}.
\begin{figure}[!h]  
\begin{center}  
\psfrag{A}{{\small $t_n$}}
\psfrag{B}{{\small $\col(t_n)$}}
\psfrag{D}{{\small $f_n$}}
\psfrag{E}{{\small $\lambda^{k} \spt(f_{n+k})$}}
\psfrag{F}{{\small $\omega^k(\col(f_{n+k}))$}}
\includegraphics[width=9cm]{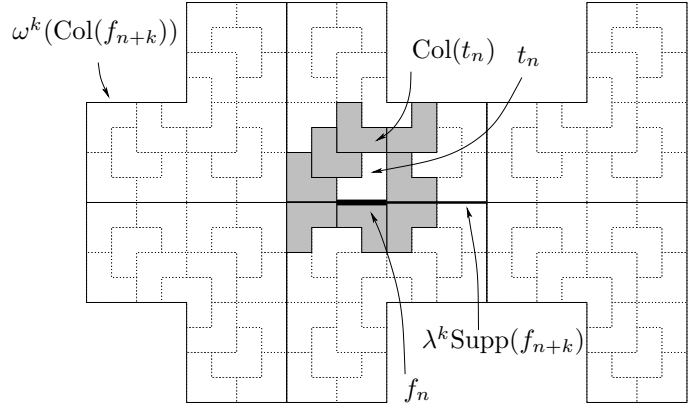}  
\end{center}
\caption{\small{Border forcing for faces.}}
\label{bratteli10.fig-facesborderforcing}
\end{figure}

As there are finitely many tiles and faces up to translation, one can choose an integer $k$ such that this property holds for any such $z$ and $x$.
\end{proof}

\begin{proposi} 
\label{bratteli10.prop-translation}
\( x\sim y\) if and only if there exists \(a(x,y) \in \RM^d\) such that \(T_x = T_y + a(x,y)\).
\end{proposi}
\begin{proof}
Consider \((x, y)\in\rel_\Bb\), and pick \(z\in  \gbo(x) \cap \gbo(y) \cap \gpath_{m,\infty}\) for some $m$ large enough for the intersection to be non empty.
And denote by $f_n$ the face corresponding to the vertex in $z$ at depth $n$.
By Proposition~\ref{bratteli10.prop-Robingpath}, $z$ defines a unique tiling $\tilde{T}_z \in \Omega$.

Let $t_n$ denote the tile corresponding to the vertex in $x$ at depth $n$.
By Lemma~\ref{bratteli10.lem-facesborderforcing}, for all $n\ge m$,  \(\col(t_n)\) appears in \(\omega^{k+n}(\col(f_{k+n}))\), see Figure~\ref{bratteli10.fig-facesborderforcing}.
This means that for all $n$, $\rob_n^c(x)$ appears in $\tilde{T}_z$: \( \rob_n^c(x) +u_n \subset \tilde{T}_z\) for some \(u_n \in \RM^d\).
But \(\rob_n(x)^\punc = \rob_m(x)^\punc\) for all \(n\ge m\), hence \(u_n = u_m, n\ge m\).
Since the sequence \((\rob_n(x)^c)_{n\ge m}\) defines $T_x$ we deduce \(\tilde{T}_z=T_x+u_m\).
The same argument shows that \(\tilde{T}=T_y+v_{m}\) for some \(v_{m} \in \RM^d\).
Therefore \(T_x = T_y + u_m - v_{m}\), and set \(a(x,y)=u_m-v_{m}\).

Conversely, consider \(T,T' \in \Xi\) with \(T=T'+a\) for some $a$.
There exists $m$ such that for all $n\ge m$,  the intersection \(\rob_n(x_T)\cap (a+\rob_n(x_{T'}))\) is non empty, and thus contains the substitute \(\omega^n(f_n)\) of a face.
This defines a sequence of faces \((f_n)_{n\ge m}\).
For $n\ge m$, $f_n$ appears on the ``common boundary'' of $t_n$ and $t'_n$ and thus \(f_n \subset \omega(f_{n+1})\).
Hence \((f_n)_{n\ge m}\) defines a generalized path $z$ in \(\gbo(x_T)\cap\gbo(x_{T'})\), and by Remark \ref{bratteli10.rem-eqrel} this proves \(x_T \sim x_{T'}\).
\end{proof}

We state two technical results that we shall need for Theorem \ref{bratteli10.thm-eqrel}.
\begin{coro} 
\label{bratteli10.cor-bda}
If \( \gbo(x) \cap \gbo(y) \cap \gpath_{m,\infty} \ne \emptyset\), then \(a(x,y) \le c \lambda^m\), where $c>0$ is a constant that does not depend on $x,y,m$.
\end{coro}
\begin{proof}
We use the notation of the first two paragraphs in the proof of Proposition \ref{bratteli10.prop-translation}.
The vector $u_m$ links the puncture of $\grob_m^c(z)$ to the puncture of a copy of $\rob_m(x)$ that it contains.
Hence \(u_m=\lambda^m u\), where $u$ is the vector linking the puncture of \(\col(f_m)\) to the puncture of a copy of $t_m$ that it contains.
The same argument shows that \(v_m = \lambda^m v\), where $v$ is the vector linking the puncture of \(\col(f_m)\) to the puncture of one of its tiles.
Now the distance from the puncture of a tile to the puncture of any of its faces is bounded above by the outer radius $R$ of the tiles, so we have \(|u|, |v|\le R\).
Hence \( |a(x,y) | = |u_m-v_m|\le |u_m| + |v_m| \le 2R\lambda^m\).
\end{proof}

\begin{lemma}
\label{bratteli10.lem-bda}
If \(T=T'+a\), then \(\gbo(x_T)\cap\gbo(x_{T'})\cap\tilde{\Pi}_{m_a,\infty} \ne \emptyset\) for some \(m_a \in \NM\) that only depends on $a$.
\end{lemma}
\begin{proof}
For each $n$, the origin and the point \(-\lambda^{-n} a\) are punctures of tiles in \(\omega^{-n}(T)\).
For all $n$ large enough those tiles must then intersect, so they have at least a common face $f_n$.
The distance between the punctures of two neighboring tiles is bounded below by $2r$, where $r$ is the inner radius of the tiles.
Let $m_a$ be the smallest integer $n$ for which \(\lambda^{-n} a \le 2r\).
The sequence of faces \((f_n)_{n\ge m_a}\) defines a generalized path, and this completes the proof.
\end{proof}

We can now state our main theorem.
\begin{theo}
\label{bratteli10.thm-eqrel}
The two equivalence relations $\rel_\Bb$ and $\rel_\Xi$ are homeomorphic:
\[
\rel_\Bb \ \cong \ \rel_\Xi \,.
\]
The homeomorphism is induced by the Robinson map \(\rob : \rinfpath \rightarrow \Xi\) of Theorem \ref{bratteli10.thm-homeotrans}.
\end{theo}
\begin{proof}
Consider the map \( \rob^\ast : \rel_\Bb \rightarrow \rel_\Xi\), given by
\[
\rob_\ast (x,y) = \bigl( \rob(x), \rob(y) \bigr) = \bigl( T_x, T_y ) = \bigl( T_x, T_x+a(x,y) \bigr) \,.
\]
By Proposition \ref{bratteli10.prop-translation}, $\rob^\ast$ is a bijection.
We now show that $\rob_\ast$ and its inverse are continuous.

Consider a sequence \((x_n,y_n)_{n\in\NM}\) that converges to \((x,y)\) in $\rel_\Bb$, and let \(a_n = a(x_n,y_n)\).
By definition of the convergence in $\rel_\Bb$, there is an $m\in\NM$ such that  \(\gbo(x_n)\cap\gbo(y_n)\cap\tilde{\Pi}_{m,\infty} \ne \emptyset\) for $n$ large enough.
Hence by Corollary \ref{bratteli10.cor-bda}, the sequence \((a_n)_{n\in\NM}\) is bounded.
By finite local complexity, it can only take finitely many values.
Any convergent subsequence must therefore be eventually stationary; and its limit must be $a(x,y)$.
Hence the sequence \((a_n)_{n\in\NM}\) has a unique accumulation point, namely $a(x,y)$, and therefore converges to $a(x,y)$ in $\RM^d$.
Since $\rob$ is a homeomorphism, \(T_{x_n} \rightarrow T_x\) and \(T_{y_n} \rightarrow T_y\) in $\Xi$.
Hence \((T_{x_n},T_{y_n}) \rightarrow (T_x,T_y)\) in $\rel_\Xi$, and this proves that $\rob_\ast$ is continuous.

Conversely, if \(((T_n,T'_n = T_n + a_n))_{n\in\NM}\) converges to \((T,T'=T+a)\) in $\rel_\Xi$, then \(a_n = a\) for all $n$ large enough.
By Lemma \ref{bratteli10.lem-bda}, there exists $m_a\in\NM$ such that \(\gbo(x_{T_n})\cap\gbo(x_{T'_n})\cap\tilde{\Pi}_{m_a,\infty} \ne \emptyset\) for all $n$ large enough.
Since $\rob$ is a homeomorphism, \(x_{T_n}\rightarrow x_T\) and \(x_{T'_n}\rightarrow x_{T'}\) in $\rinfpath$.
And therefore \((x_{T_n},x_{T'_n}) \rightarrow (x_{T}, x_{T'})\) in $\rel_\Bb$.
This proves that $(\rob_\ast)^{-1}$ is continuous.
\end{proof}

Since $\rel_\Xi$ is an {\it \'etale} equivalence relation (Definition~\ref{bratteli10.def-etale}), we have the immediate corollary.
\begin{coro} 
\label{bratteli10.cor-etale}
The equivalence relation $\rel_\Bb$ is {\em \'etale}.
\end{coro}
As a consequence of Lemma~\ref{bratteli10.lem-facesborderforcing} and Corollary~\ref{bratteli10.cor-bda} in particular, one can check that a base of $\rel_\Bb$-sets for the {\it \'etale} topology of $\rel_\Bb$ is given by the following
\begin{multline}
\label{bratteli10.eq-rsetsB}
O_{\gamma\gamma', \eta, m} = 
\Bigl\{
(x,y) \in \rel_\Bb \ : \ x\in[\gamma], \; y\in[\gamma'],\; \\
\gbo(x)\cap\gbo(y)\cap\tilde{\Pi}_{m-k,\infty}  \neq \emptyset \,, 
\gbo(x)\cap\gbo(y)\cap\tilde{\Pi}_{m-k,m}  = \eta
\Bigr\}
\end{multline}
where \(\gamma, \gamma'\in \Pi_{0, m}^d\), \(m > k\), with $k$ the parameter of border forcing for the faces as defined in Lemma~\ref{bratteli10.lem-facesborderforcing}, and where $\eta$ is a generalized path of length $k$ in \(\gbo(\gamma_{[m-k,k]})\cap\gbo(\gamma'_{[m-k,k]})\).
One can check that those sets are compatible with the topology of $\rel_\Bb$.
Also, one sees that the restrictions of the range and source maps to those sets are homeomorphic as follows.

\begin{figure}[!h]  
\begin{center}  
\psfrag{A}{{\small $\rob_{m-k}(\gamma_{[0,m-k]})$}}
\psfrag{B}{{\small $\rob_{m-k}(\gamma'_{[0,m-k]})$}}
\psfrag{C}{{\small $\rob^c_{m-k}(\gamma_{[0,m-k]})$}}
\psfrag{D}{{\small $\rob^c_{m-k}(\gamma'_{[0,m-k]})$}}
\psfrag{I}{{\small $\rob_m(\gamma)$}}
\psfrag{J}{{\small $\rob_m(\gamma')$}}
\psfrag{G}{{\small $\grob^c_m(\eta)$}}
\psfrag{E}{{\small $\lambda^{m-k} \spt(z_{m})$}}
\psfrag{F}{{\small $\lambda^m \spt(z_{m-k})$}}
\includegraphics[width=10cm]{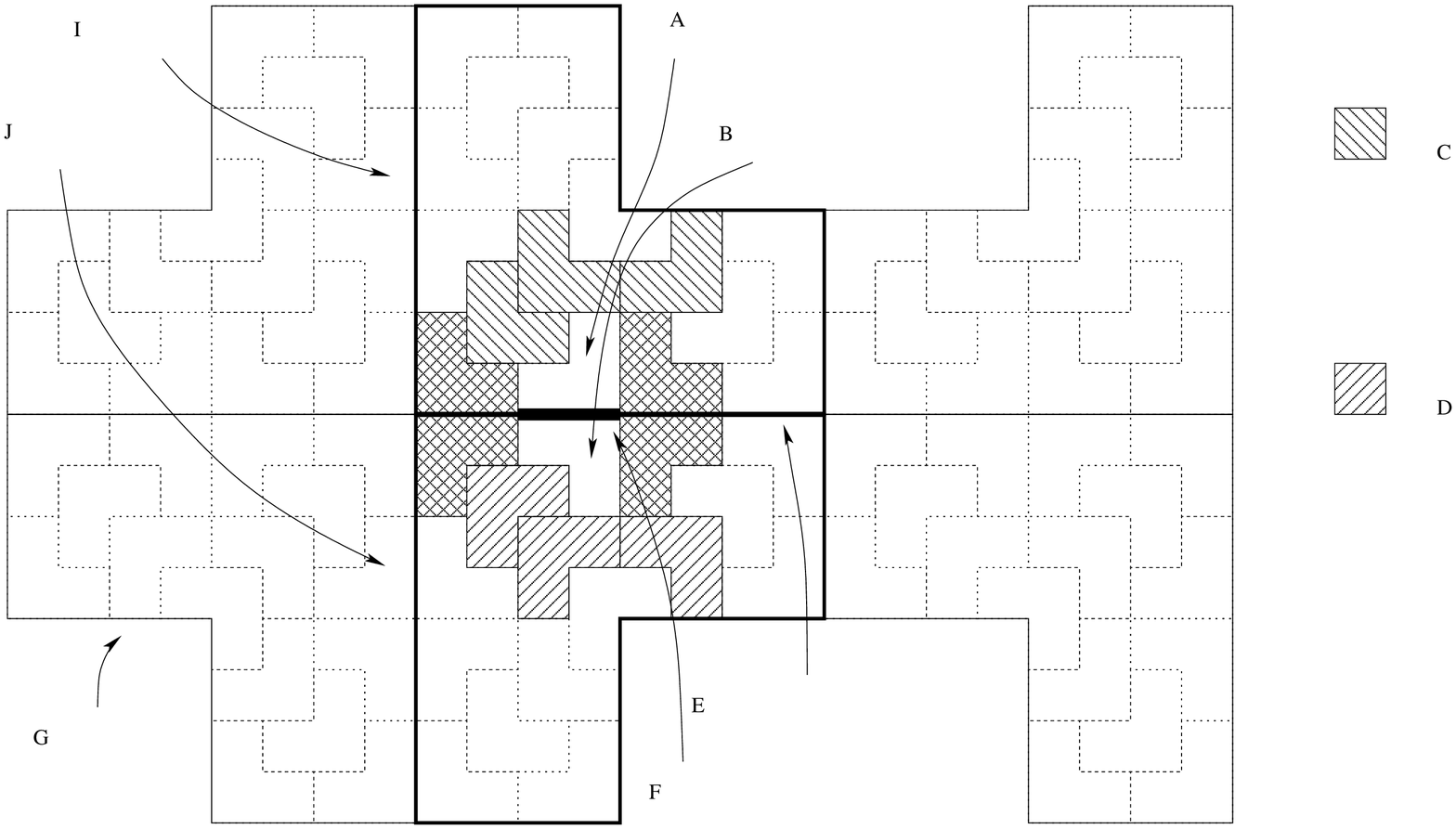}  
\end{center}
\caption{\small{Geometric condition for the $\rel_\Bb$-sets.}}
\label{bratteli10.fig-rsets}
\end{figure}

The condition in equation~\eqref{bratteli10.eq-rsetsB} says that the first $k+1$ vertices of any generalized path in \(\gbo(x)\cap\gbo(y)\cap\tilde{\Pi}_{m-k,\infty}\) are determined by $\eta=(z_{m-k}, z_{m-k+1}, \ldots, z_m)$.
Hence for any \((x,y) \in O_{\gamma\gamma', \eta, m}\), the patches \(\phi_{m-k}^c(x) = \phi_{m-k}^c(\gamma)\) and \(\phi_{m-k}^c(y) = \phi_{m-k}^c(\gamma')\) appear in $\grob_m^c(\eta)$, hence in both \(\phi_{m}(x) = \phi_{m}(\gamma)\) and \(\phi_{m}(y) = \phi_{m}(\gamma')\), and at respective positions which are uniquely determined by $\eta$.
Therefore the distance between the tilings $T_x$ and $T_y$ is uniquely determined by the paths $\gamma$, $\gamma'$ and the generalized path $\eta$.
In other words, given any compatible $x\in [\gamma]$, there exists a unique $y\in[\gamma']$ such that \((x,y) \in O_{\gamma\gamma',\eta, m}\).
See Figure~\ref{bratteli10.fig-rsets} for an illustration.
Notice that for fixed $\gamma, \gamma', m$, two sets \(O_{\gamma\gamma', \eta_1, m}\) and \(O_{\gamma\gamma', \eta_2, m}\) are either equal or disjoint (this is because $\tilde{\phi}^c$ is not one-to-one: two generalized paths $\eta_1$ and $\eta_2$ might determine the same faces on the boundary of \(\phi_{m}(\gamma)\) and \(\phi_{m}(\gamma')\)).

Another immediate corollary of Theorem~\ref{bratteli10.thm-eqrel} is the following.
\begin{coro} 
\label{bratteli10.cor-groupoid}
The groupoid of the equivalence relation $\rel_\Bb$ is homeomorphic to $\Gamma_\Xi$.
\end{coro}

\end{document}

%% file: macros.tex
\usepackage[T1]{fontenc}
\usepackage[utf8]{inputenc}
\usepackage[cm]{aeguill}
\usepackage[english]{babel}

\usepackage[all]{xy}
\usepackage{amsmath}
\usepackage{amssymb}
\usepackage{amsthm}
\usepackage{enumerate}

\usepackage{graphicx}
\usepackage{psfrag}

\usepackage{geometry}

\newtheorem{theo}{Theorem}[section]
\newtheorem{defini}[theo]{Definition}
\newtheorem{proposi}[theo]{Proposition}
\newtheorem{lemma}[theo]{Lemma}
\newtheorem{coro}[theo]{Corollary}
\newtheorem{rem}[theo]{Remark}

\newtheorem{hypo}[theo]{Hypothesis}

\newtheorem{notat}[theo]{Notation}

\newcommand{\Aa}{{\mathcal A}}
\newcommand{\Bb}{{\mathcal B}}

\newcommand{\Ee}{{\mathcal E}}
\newcommand{\Ff}{{\mathcal F}}

\newcommand{\Hh}{{\mathcal H}}

\newcommand{\Rr}{{\mathcal R}}

\newcommand{\Vv}{{\mathcal V}}

\newcommand{\NM}{{\mathbb N}}

\newcommand{\RM}{{\mathbb R}}

\newcommand{\ZM}{{\mathbb Z}}

\newcommand{\Css}{$C^{\ast}$-algebras }            

\newcommand{\eaf}{\stackrel{\mbox{\tiny\sc AF}}{\sim}}
\newcommand{\etail}{\stackrel{\mbox{\tiny\rm tail}}{\sim}}
\newcommand{\col}{\text{\rm Col}}                    
\newcommand{\punc}{\text{\rm punc}}                  
\newcommand{\spt}{\text{\rm Supp}}                   
\newcommand{\bd}{\text{\rm bdim}}                    
\newcommand{\dist}{\text{\rm dist}}                  
\newcommand{\vtx}{\Vv}                  
\newcommand{\edg}{\Ee}                  
\newcommand{\rel}{\Rr}                  
\renewcommand{\tilde}{\widetilde}                  
\newcommand{\infpath}{\Pi_{\bullet,\infty}}                  
\newcommand{\rinfpath}{\Pi_{\infty}}                  
\newcommand{\gpath}{\tilde{\Pi}}                  
\newcommand{\infgpath}{\tilde{\Pi}_{\bullet,\infty}}                  
\newcommand{\bo}{\text{\rm b}}                  
\newcommand{\gbo}{\tilde{\text{\rm b}}}                  

\renewcommand{\phi}{\varphi}                  
\newcommand{\rob}{\phi}                               
\newcommand{\grob}{\tilde{\phi}}                  

\renewcommand{\tilde}{\widetilde}

\newcommand{\face}[1]{\proto^{#1}}
\newcommand{\proto}{\Aa}
\newcommand{\vois}[1]{[#1]}
\newcommand{\mbrat}[1]{\Bb^{#1}}   
\newcommand{\rac}{\circ}
\newcommand{\ra}{\rightarrow}
\newcommand{\lra}{\longrightarrow}

\newcommand{\restr}[2]{#1_{\vert #2}}

\newcommand{\escaping}{\mathcal{S}}
\newcommand{\hor}[2][]{\mathcal{H}_{#2}^{#1}}  
\newcommand{\chemin}[2]{\Pi_{#1}^{#2}}

\newcommand{\nr}[1]{\left\Vert #1 \right\Vert}